\documentclass[10pt,a4paper,twoside]{scrartcl}
\usepackage[utf8]{inputenc}
\usepackage[english]{babel}
\usepackage{lmodern}
\usepackage[left=2cm,right=2cm,top=2.5cm,bottom=2.5cm]{geometry}
\usepackage{amsmath,amssymb,mathtools,amsfonts}
\usepackage{graphicx}
\usepackage{subcaption}
\usepackage[]{algorithm2e}
\usepackage{paralist}

\usepackage{hyperref}
\usepackage[capitalize, nameinlink, noabbrev]{cleveref}
\crefname{equation}{}{}
\usepackage{autonum}

\usepackage{scrlayer-scrpage}
\ohead{N. Schneider and V. Michel: A dictionary learning add-on for spherical downward continuation}

\newcommand{\Hs}[1]{\mathcal{H}_#1(\Omega)}
\newcommand{\Lp}[1]{\mathrm{L}^#1(\Omega)}
\newcommand{\T}{\mathcal{T}}
\newcommand{\ball}{\mathbb{B}}
\newcommand{\real}{\mathbb{R}}
\newcommand{\nat}{\mathbb{N}}
\newcommand{\dic}{\mathcal{D}}
\newcommand{\dicinf}{\mathcal{D}^{\mathrm{Inf}}}
\newcommand{\Can}{\mathcal{C}}
\newcommand{\RFMP}{\mathrm{RFMP}}
\newcommand{\ROFMP}{\mathrm{ROFMP}}

\newcommand{\IPMP}{\mathrm{IPMP}}
\newcommand{\N}{\mathcal{N}}

\newcommand{\SH}{\mathrm{SH}}
\newcommand{\APK}{\mathrm{APK}}
\newcommand{\APW}{\mathrm{APW}}
\newcommand{\SL}{\mathrm{SL}}

\newcommand{\SO}[1]{\mathrm{SO}(#1)}

\newcommand{\lon}{\varphi}
\newcommand{\lat}{\theta}

\newcommand{\ethree}{\varepsilon^3}

\newcommand{\trans}{\mathrm{T}}
\newcommand{\proj}[1]{\mathcal{P}_{#1}}

\title{A dictionary learning add-on for spherical downward continuation}
\author{N. Schneider$^\ast$ and V. Michel\footnote{Geomathematics Group Siegen, University of Siegen, corresponding author: naomi.schneider@mathematik.uni-siegen.de}}
\date{}

\usepackage{natbib}
\begin{document}
\maketitle

\begin{abstract}
We propose a novel dictionary learning add-on for existing approximation algorithms for spherical inverse problems such as the downward continuation of the gravitational potential.
The Inverse Problem Matching Pursuit (IPMP) algorithms iteratively minimize the Tikhonov functional in order to construct a weighted linear combination of so-called dictionary elements as a regularized approximation. A dictionary is a set that contains trial functions such as spherical harmonics (SHs), Slepian functions (SLs) as well as radial basis functions (RBFs) and wavelets (RBWs). Previously, the IPMP algorithms worked with finite dictionaries which are vulnerable regarding a possible biasing of the outcome. 
Here, we propose an additional learning technique that allows us to work with infinitely many trial functions and provides us with a learnt dictionary for future use in the IPMP algorithms. We explain the general mechanism and provide numerical results that prove its applicability and efficiency.
\end{abstract}

\paragraph{Keywords}
dictionary learning, inverse problems, gravitational potential, matching pursuits, numerical modelling, satellite geodesy

\paragraph{Acknowledgments}
The authors gratefully acknowledge the financial support by the German Research Foundation (DFG; Deutsche Forschungsgemeinschaft), projects MI 655/7-2 and MI 655/14-1. Further, we thank Dr. Roger Telschow for handing us the irregularly distributed grid.

\section{Introduction}
\label{sect:intro}
The climate is changing. In particular, the repercussions caused by the loss of water in the soil or the ice sheets are dangerous to human kind as they fuel droughts as well as the sea level rise. Thus, the geosciences are taking care to monitor the mass transport of the Earth, see e.\,g.\ \cite{Fischeretal2013-2,Flechtneretal2020,IPCC2014,Linetal2018,NasaConsensu,Saemianetal2019,Wieseetal2020} as well as the results of the DFG SPP 1257 (2006-2014) coordinated by Ilk and Kusche, see e.\,g. \cite{Kuscheetal2012}. The mass transport is obtained by modelling the gravitational potential from time-dependent satellite data, e.\,g. from the GRACE and GRACE-FO satellite mission, see e.\,g.  \cite{Devaraju2017,Flechtneretal2014,Flechtneretal2014-2,GRACEdata,Schmidtetal2008,Tapleyetal2004,GRACEdata2}.\\
We approximate the Earth's surface with the unit sphere $\Omega$.
The gravitational potential $f$ on the Earth's surface is usually expanded in SHs $Y_{n,j},\ n\in\nat_0,\ j=-n,...,n,$ such that we have 
	\begin{align}
		f = \sum_{n=0}^\infty \sum_{j=-n}^n \left\langle f,Y_{n,j}\right\rangle_{\Lp{2}}Y_{n,j}
		\label{eq:gravpotatearth}
	\end{align}	 
which holds in the $\Lp{2}$-sense. Correspondingly, on a satellite orbit $\sigma>1$, we upward continue this representation at every point $\eta \in \Omega$ via the operator $\T$ and obtain the gravitational potential $V$ outside the Earth, see e.\,g. \cite{Baur2014,Freedenetal2004,Moritz2010,Telschow2014},
	\begin{align}
		V\left(\sigma\eta\right) = \left( \T f \right) \left(\sigma\eta\right) = \sum_{n=0}^\infty \sum_{j=-n}^n \left\langle f,Y_{n,j}\right\rangle_{\Lp{2}}\sigma^{-n-1}Y_{n,j}\left(\eta\right),
		\label{eq:gravpotatorbit}
	\end{align} 
which actually holds pointwise. Naturally, we are much more interested in the inverse problem of the downward continuation, i.\,e. if we have given values of the potential $V$ at $\sigma\eta$, we are interested in the values of $f$ at $\eta$. Due to the exponentially decreasing singular values of $\T$, its inverse is not continuously dependent on the data, see e.\,g. \cite{Engletal1996,Louis1989,Michel2005,Michel2020,Rieder2003,Telschow2014}. This means, mathematically, the downward continuation is an exponentially ill-posed inverse problem and needs sophisticated methods to be tackled.\\
One possible approach are
 the IPMP algorithms, i.\,e. here the Regularized (Orthogonal) Functional Matching Pursuit (R(O)FMP) algorithm, see e.\,g. \cite{Berkeletal2011,Fischer2011,Fischeretal2012,Fischeretal2013-1,Fischeretal2013-2,Guttingetal2017,Kontak2018,Kontaketal2018-2,Kontaketal2018-1,Leweke2018,Michel2015-2,Micheletal2017-1,Micheletal2018-1,Micheletal2014,Micheletal2016-1,Schneider2020,Telschow2014}. Using discrete values of $V$, these methods iteratively build an approximation of $f$ as a best basis expansion drawn from dictionary elements. A dictionary $\dic$ is a set of trial functions. For spherical tasks like the downward continuation, it may contain SHs, SLs, RBFs and/or RBWs. Then $\dic$ contains global functions as well as localized ones. In an iteration $N$, the next dictionary element is chosen such that it minimizes the Tikhonov functional. In the ROFMP algorithm, it simultaneously fulfils to a certain extent orthogonality relations with the previously chosen basis elements as well. Further, the coefficients are updated regularly to maintain their optimality. However, the IPMP algorithms usually work with a finite, a-priori chosen subset of the infinitely many possible trial functions as its dictionary. We immediately recognize that this may bias the obtained approximation.\\
Thus, we further developed the IPMP algorithms by adding a novel dictionary learning technique. Our task at hand differs from established dictionary learning challenges because of the considered inverse problem as well as our focus on a non-discrete result. Thus, we cannot straightforwardly transfer previous dictionary learning strategies, see e.\,g. \cite{Aharonetal2006,Brucksteinetal2009,Enganetal1999,Enganetal1999-1,Pruente2008,Rubinsteinetal2010}. Nonetheless, our approach can successfully be discussed in the light of similar learning criteria (such as flexibility and a well-defined goal) and a similar theoretical starting point (i.\,e. solving a doubled minimization problem). For further information of these relations, the interested reader is referred to \cite{Schneider2020}.\\
Here, the dictionary learning add-on for the SLs, RBFs and RBWs consists of a 2-step optimization process in order to compute one optimized candidate for each type of these trial function. The process solves non-linear constrained optimization problems globally and locally. We utilize the NLOpt library for this. A candidate for the SHs can be obtained by comparing the values of the objective function for different SHs. Then we learn particular functions as well as a maximal SH degree. All of these candidates together constitute again a finite dictionary and we can proceed with the usual routine of the IPMP algorithms. After termination, we obtain an approximation of $f$ in an optimized best basis whose elements can be re-used as a learnt dictionary in future runs of the IPMP algorithms. The IPMP algorithms which include the novel learning add-on are called the Learning Inverse Problem Matching Pursuit (LIPMP) algorithms, i.\,e. the Learning Regularized (Orthogonal) Functional Matching Pursuit (LR(O)FMP) algorithm. Note that the LRFMP algorithm here is an advanced version of the one presented in \cite{Micheletal2018-1}.\\
In the sequel, we formally but shortly introduce the SHs, SLs, RBFs and RBWs as well as a dictionary in \cref{sect:basics}. Then we explain the idea of the LIPMP algorithms which includes an overview of the IPMP algorithms and the novel learning add-on in \cref{sect:LIPMPs}. Particularly, we focus on how the 2-step optimization process fits into the routine of the IPMP algorithms. We also summarize shortly some theoretical aspects of the LIPMP algorithms. Then we show the applicability and efficiency of the add-on as well as the learnt dictionary in a series of experiments in \cref{sect:numresults}.\\
This paper is based on \cite{Schneider2020} and was in principle presented at the EGU2020: Sharing Geoscience Online \cite{Schneider2020EGU}. 

\section{Basics}
\label{sect:basics}
Let $\real$ be the set of all real numbers, $\real^d$ be the real, $d$-dimensional vector space, $\nat$ be the set of all positive integers and $\nat_0$ be that of all non-negative integers. We denote $\Omega\coloneqq \{x\ \in \real^3 :\ |x|=1\}$ as the unit sphere and $\ball \coloneqq \ball^3 \coloneqq \{x\ \in \real^3 :\ |x|< 1\}$ the open unit ball. Furthermore, we can represent $\eta(\lon,t) \in \Omega$ as usual, see e.\,g. \cite{Michel2013}, via the longitude $\lon \in [0,2\pi[$ and the latitude $t = \cos(\lat) \in [-1,1]$ for $\lat \in [0,\pi]$.

\subsection{Spherical Harmonics}
Spherical harmonics (SHs) $Y_{n,j}$ are (global) polynomials on $\Omega$. They have a distinct degree $n \in \nat_0$ and order $j=-n,...,n$. In practice, we usually choose the fully normalized SHs 
\newpage
	\begin{align}
		Y_{n,j}(\eta(\lon,t)) \coloneqq p_{n,j}P_{n,|j|}(t) \left\{ \begin{matrix} \cos(j\lon),&j\leq 0\\ \sin(j\lon),&j> 0\end{matrix} \right.
		\label{def:fnsh}
	\end{align}
for $\eta(\lon,t) \in \Omega$, an $\Lp{2}$-normalization $p_{n,j}$ and the associated Legendre functions $P_{n,|j|}$. An example is given in \cref{fig:TF1}. For further details, the reader is referred to, e.\,g., \cite{Freedenetal1998,Freedenetal2013-1,Freedenetal2004,Freedenetal2009,Michel2013,Mueller1966}.

\subsection{Slepian Functions}
Slepian functions (SLs) are band-limited, spatially optimally localized trial functions. Here, the spatial localization region shall be a spherical cap. A spherical cap can be parametrized by the parameter $c = \cos(\lat) \in [-1,1],$ where $\lat$ is the polar angle between vectors pointing at the apex and at an arbitrary point of the base, and its centre $A(\alpha,\beta,\gamma)\ethree \in \Omega$, where $\alpha,\ \gamma \in [0,2\pi]$ and $\beta \in [0,\pi]$ denote the Euler angles, $A\in\SO{3}$ is a rotation matrix and $\ethree =(0,0,1)^\trans$ is the North Pole. Then, for $k=1,...,(L+1)^2,$ we set 
	\begin{align}
		g^{(k,L)}(R,\eta) \coloneqq \sum_{l=0}^L\sum_{m=-l}^l g^{(k,L)}_{l,m}(R) Y_{l,m}(\eta),\ \eta \in \Omega,
		\label{def:sl}
	\end{align}
where $R \coloneqq R(c,A(\alpha,\beta,\gamma)\ethree) \in \overline{\ball}^4 \coloneqq \left\{x \in \real^4\ \colon\ |x|\leq 1\right\}$ stands for the localization region and $L \in \nat_0$ is the band-limit. The Fourier coefficients $g^{(k,L)}_{l,m}(R),\ l=0,...,L,\ m=-l,...,l,$ are obtained from the eigenvectors of the related algebraic eigenvalue problem of optimizing a band-limited, in SHs expanded function in the region $R$. The superscript $k$ enumerates the $(L+1)^2$ Slepian functions with band-limit $L$. An example is given in \cref{fig:TF1}. Note that a commuting operator provides a stable computation of these values if the localization region is a spherical cap. For further details, the reader is referred to, e.\,g., \cite{Albertellaetal1999,Gruenbaumetal1982,Michel2013,Seibert2018,Simonsetal2006-1}.
	
\subsection{Radial Basis Functions and Wavelets}
As examples for localized trial functions, we consider Abel--Poisson kernels $K(x,\cdot)$ (APKs) and wavelets $W(x,\cdot)$ (APWs) due to their closed form. That means, we have
	\begin{align}
		K(x,\eta) &\coloneqq \frac{1-|x|^2}{4\pi(1+|x|^2-2x\cdot\eta)^{3/2}}
		\label{def:apk}
	\intertext{and}
		W(x,\eta) &\coloneqq K(x,\eta) - K(|x|x,\eta)
		\label{def:apw}
	\end{align}
for $\eta \in \Omega$ and with the characteristic parameter $x\in\ball$. The kernels act as low pass filters whereas the wavelets are band pass filters. Examples are given in \cref{fig:TF2}. For further details, the reader is referred to, e.\,g., \cite{Freedenetal1998,Freedenetal2004,Freedenetal1998-1,Freedenetal2009,Freedenetal1996,Michel2013,Windheuser1995}.\\

\subsection{A dictionary}
A dictionary is a set of trial functions. We first define subsets for each type of trial function under investigation. Here, we set 
	\begin{align}
		[N]_{\SH} &\coloneqq \left\{ (n,j) \in N \subseteq \N \right\},
		\label{def:tfc:sh}\\
		\N &\coloneqq \left\{(n,j)\ \colon\ n\in \nat_0,\ j=-n,...,n\right\}
	\intertext{for SHs,} 
		[S]_{\SL} &\coloneqq \left\{ g^{(k,L)}(R,\cdot)\ \colon\ (R,(k,L)) \in \overline{\ball}^4 \times \mathcal{L}\right\}
		\label{def:tfc:sl}\\
		\mathcal{L} &\coloneqq \left\{(k,L)\ \colon\ L\in \nat_0,\ k=1,...,(L+1)^2\right\}
	\intertext{for SLs,} 
		[B_K]_{\APK} &\coloneqq \left\{ x \in \ball\ \colon\ x \in B_K \subseteq \ball \right\}
		\label{def:tfc:apk}
	\intertext{for APKs and} 
		[B_W]_{\APW} &\coloneqq \left\{ x \in \ball\ \colon\  x \in B_W \subseteq \ball \right\}
		\label{def:tfc:apw}
	\end{align}
for APWs. Such subsets are called trial function classes. The union of the defined trial function classes gives us the dictionary $\dic$:
	\begin{align}
		\dic \coloneqq [N]_{\SH} \cup [S]_{\SL} \cup [B_K]_{\APK} \cup [B_W]_{\APW},
		\label{def:dic}
	\end{align}
confer \cite{Michel2020,Micheletal2018-1,Schneider2020}.
In general, it is not necessary that $\dic$ is finite. However, we emphasize an infinite dictionary as $\dicinf$. Note that, depending on the actual choice of $[N]_\SH,\ [B_K]_\APK$ and $[B_W]_\APW$, $\dic$ may be complete in certain function spaces like $\Lp{2}$, see e.\,g. \cite{Freedenetal1998,Michel2013,Schneider2020} and the previously mentioned references regarding the SHs, APKs and APWs. This holds in particular for $\dicinf$.

\section{The LIPMP algorithms}
\label{sect:LIPMPs}
We consider the downward continuation of the gravitational potential from satellite data to the Earth's surface. That is, mathematically speaking, we consider the ill-posed inverse problem $y=\T_\daleth f$ with the data $y \in \real^\ell,\ y_i = V(\sigma\eta^i)$, the satellite height $\sigma >1$, grid points $\eta^i \in \Omega$ for $i=1,..., \ell$ at the Earth's surface and the operator $\T_\daleth f \coloneqq ((\T f)(\sigma\eta^i))_{i=1,...,\ell}$ with the upward continuation operator $\T$ as in \cref{eq:gravpotatorbit}. Thus, $\T_\daleth$ is the corresponding evaluation operator of the upward continuation operator for a $\ell$-dimensional discretized grid. Our task is to approximate the gravitational potential $f$ at the Earth's surface $\Omega$. \\
The LIPMP algorithms introduce an add-on to the established IPMP algorithms. The remaining routines coincide. Thus, we give a short overview of the strategy.

\subsection{The underlying IPMP algorithms}
Due to the ill-posedness of the inverse problem, we need to consider a regularization for the downward continuation. The IPMP algorithms utilize a Tikhonov regularization. In particular, using $f_0 \equiv 0$, the methods iteratively build a linear combination of weighted dictionary elements $d_n \in \dic$:
	\begin{align}
		f_N &= \sum_{n=1}^N \alpha_n d_n
		\label{eq:fn}
	\intertext{in the case of the RFMP algorithm and}
		f_N^{(N)} &= \sum_{n=1}^N \alpha_n^{(N)} d_n,
		\label{eq:fnN}
	\end{align}
in the case of the ROFMP algorithm. The superscript $(N)$ refers to the update of the coefficients in each iteration step due to the orthogonality process. Note that, in practice, we usually consider the iterated (L)ROFMP algorithm which restarts the orthogonality technique after a prescribed number of iterations for practical as well as theoretical reasons. For readability, here, we do without the additional subscripts that are necessary for the restart process. \\
The respective residuals are
	\begin{align}
	R^{N+1} &\coloneqq R^N - \alpha_{N+1}\T_\daleth d_{N+1}
		\label{def:RN}
	\intertext{for the RFMP algorithm and, in the case of the ROFMP algorithm,}
	R^{N+1} &\coloneqq R^N - \alpha_{N+1}\proj{\mathcal{V}_N^\perp} \T_\daleth d_{N+1},
		\label{def:RNO}
	\intertext{where $\proj{\mathcal{V}_N^\perp}$ is the orthogonal projection onto}
	\mathcal{V}_N^\perp &\coloneqq \left(\mathrm{span} \{\T_\daleth d_n\ \colon n=1,...,N\}\right)^{\perp_{\real^\ell}}.
		\label{def:VNO}
	\end{align}
In both cases, we have $R^0 = y-\T_\daleth f_0 = y$. In each iteration $N$, the weights $\alpha_{N+1} \in \real$ and $\alpha_{N+1}^{(N+1)} \in \real$, respectively, as well as the basis function $d_{N+1} \in \dic$ are chosen such that, for $\lambda>0$, the Tikhonov functional 
	\begin{align}
	\left\| R^N - \alpha \T_\daleth d \right\|^2_{\real^\ell} + \lambda\left\|f_N+\alpha d\right\|^2_{\Hs{2}}
		\label{eq:TFNO}
	\end{align}
for the RFMP algorithm and, for the ROFMP algorithm,
	\begin{align}
	\left\| R^N - \alpha\proj{\mathcal{V}_N^\perp}\T_\daleth d \right\|^2_{\real^\ell} + \lambda\left\|f_N^{(N)}+\alpha \left(d-b_n^{(N)}(d) \right)\right\|^2_{\Hs{2}},
		\label{eq:TFO}
	\end{align}
where
	\begin{align}
	b_n^{(N)}(d) \coloneqq \sum_{n=1}^N \beta_n^{(N)}(d) d_n
		\label{def:bnN}
	\end{align}
with the projection coefficients $\beta_n^{(N)}(d)$, is minimized. The projection coefficients are given by 
\begin{align}
				\beta_N^{(N)}(d) \coloneqq \frac{\left\langle T_\daleth d, \proj{\mathcal{V}^\trans_{N-1}} T_\daleth d_N \right\rangle_{\real^\ell}}{\left\| \proj{\mathcal{V}^\trans_{N-1}} T_\daleth d_N\right\|_{\real^\ell}^2} 
				\qquad \textrm{and} \qquad
				\beta_n^{(N)}(d) \coloneqq \beta_n^{(N-1)}(d) - \beta_N^{(N)}(d) \beta_n^{(N-1)}(d_n)
				\label{eq:betan_1}
			\end{align}
for $n=1,...,N-1.$ For the penalty term, we use the norm of the Sobolev space $\Hs{2} \subset \Lp{2}$ which is the completion of the set of square-integrable functions whose $(n+0.5)^4$-weighted $\Lp{2}$-norm is finite, see e.\,g. \cite{Freedenetal1998,Michel2013}. \\
In practice, we determine the basis element $d_{N+1}$ as the maximizer of 
	\begin{align}
		\RFMP(d;N) &\coloneqq \frac{\left( \left\langle R^N, \T_\daleth d\right\rangle_{\real^\ell} - \lambda\left\langle f_N,d \right\rangle_{\Hs{2}} \right)^2}{\left\|\T_\daleth d\right\|_{\real^\ell}^2 + \lambda\left\|d\right\|^2_{\Hs{2}}}
		\eqqcolon \frac{\left(A_N(d)\right)^2}{B_N(d)}
		\label{def:RFMP(d;N)}
	\intertext{and}
		\ROFMP(d;N) &\coloneqq \frac{\left( \left\langle R^N, \proj{\mathcal{V}_N^\perp}\T_\daleth d \right\rangle_{\real^\ell} - \lambda\left\langle f_N^{(N)},d-b_n^{(N)}(d) \right\rangle_{\Hs{2}} \right)^2}{\left\| \proj{\mathcal{V}_N^\perp}\T_\daleth d\right\|_{\real^\ell}^2 + \lambda\left\|d-b_n^{(N)}(d)\right\|^2_{\Hs{2}}}
		\eqqcolon \frac{\left(A_N^{(N)}(d)\right)^2}{B_N^{(N)}(d)},
		\label{def:ROFMP(d;N)}
	\end{align}
respectively. Then we obtain the weights via 
	\begin{align}
	\alpha_{N+1} &\coloneqq \frac{A_N(d_{N+1})}{B_N(d_{N+1})}
		\label{def:alphan}
	\intertext{and}
	\alpha_{N+1}^{(N+1)} &\coloneqq \frac{A_N^{(N)}(d_{N+1})}{B_N^{(N)}(d_{N+1})},
		\label{def:alphanN}
	\end{align}
respectively. The IPMP algorithms most commonly terminate if the relative data error falls below a certain threshold like the noise level or if a certain number of iterations is reached. For more details on the IPMP algorithms, the reader is referred to \cite{Berkeletal2011,Fischer2011,Fischeretal2012,Fischeretal2013-1,Fischeretal2013-2,Guttingetal2017,Kontak2018,Kontaketal2018-2,Kontaketal2018-1,Leweke2018,Michel2015-2,Micheletal2017-1,Micheletal2018-1,Micheletal2014,Micheletal2016-1,Schneider2020,Telschow2014}.\\
The dictionary $\dic$ is finite in most of the previous publications on an IPMP algorithm as the use of an infinite dictionary $\dicinf$ was an open question at that time. Then the maximizer of \cref{def:RFMP(d;N)} and \cref{def:ROFMP(d;N)}, respectively, is found by evaluating the objective function for all dictionary elements and choosing the maximum amongst them. However, a finite dictionary may bias the approximation $f_N$ and $f_N^{(N)}$, respectively. To enable the use of an infinite dictionary and provide an automation of selecting a finite dictionary (i.\,e. learn a dictionary), the LIPMP algorithms were developed.

\subsection{The learning add-on}
We first consider how an infinite dictionary $\dicinf$ can be introduced into the routine of the IPMP algorithms. From this, the learnt dictionary follows naturally.\\
The infinite dictionary $\dicinf$ is defined by
	\begin{align}
		\dicinf \coloneqq [\widetilde{N}]_{\SH} \cup [S]_{\SL} \cup [B_K]_{\APK} \cup [B_W]_{\APW}
		\label{def:dicinf}
	\end{align}
	with
	\begin{align}
		\widetilde{N} \coloneqq \{(n,j)\ \colon\ n \leq \overline{N},\ j=-n,...,n\},\qquad
		S \coloneqq  \ball^4 \times \mathcal{L} \qquad \textrm{and}\qquad
		B_K \coloneqq 	B_W \coloneqq \ball
	\end{align}
for fixed $\overline{N},\ L \in \nat_0$. The trial function class of the SH is still finite. This is due to the discrete nature of the characteristic parameters, the degree $n$ and the order $j$. Nonetheless, the other trial function classes are truly infinite. This means that the choice of the parameters of the centre $\alpha,\ \beta,\ \gamma$ and the size $c$ of the spherical cap for a Slepian function are arbitrary, while their band-limit is fixed and finite in analogy to the SH choice. Moreover, also the centres of the radial basis functions and wavelets can be chosen from all points in the ball $\ball$.\\
The main obstacle for using $\dicinf$ is the determination of the maximizer of \cref{def:RFMP(d;N)} and \cref{def:ROFMP(d;N)}, respectively, in the truly infinite trial function classes. For this, we introduce an additional optimization step into the routine. In particular, in this step, we determine a finite dictionary of (optimized) candidates $\mathcal{C}$ from the infinite $\dicinf$. The set $\mathcal{C}$ contains one optimum per function class and we can then proceed as usual. Therefore, after termination (which correspondingly obeys to the same rules as in the IPMP algorithms), we obtain an approximation $f_N$ and $f_N^{(N)}$, respectively, in a best basis of optimized dictionary elements. The latter constitute the learnt (finite) dictionary which can be used in future runs of the IPMP algorithms.\\
Due to the different nature of the classes, we distinguish a strategy for the SH and the remaining three trial function classes. The approach to learn SHs was already explained to a certain extent in \cite{Micheletal2018-1}. For completeness and some additional insights, we summarize it here again. We propose to learn a maximal SH degree as well as particular SHs simultaneously. The trial function class $[\widetilde{N}]_\SH$ includes all SHs up to a degree $\overline{N}$ (see \cref{def:dicinf}). The value of $\overline{N}$ should be chosen much larger than it is sensible for the data $y$. Then, we again follow the previous approach and determine the values of \cref{def:RFMP(d;N)} and \cref{def:ROFMP(d;N)}, respectively, for all SHs up to degree $\overline{N}$. Hence, after termination, we have a certain set of SHs that are used in the representation \cref{eq:fn} and \cref{eq:fnN}, respectively. Firstly, most likely, the algorithms will have determined a smaller, properly learnt maximal SH degree $\nu \in \nat,\ \nu < \overline{N},$ on its own. Secondly, we have a set of distinct SHs used in the approximation and, thus, contained in the learnt dictionary. Note that this approach demands a finite starting dictionary $\dic^\mathrm{s}$ in the LIPMP algorithms which contains (at least) $[\widetilde{N}]_\SH$.\\
For the remaining trial function classes of the SLs, APKs and APWs, we determine each candidate by solving non-linear constraint optimization problems. Note that, we use $N$ for the iteration number next. The objective functions in the $N$-th iteration are $\IPMP(d(z);N)$ where $d(z)$ denotes a SL, APK or APW, respectively, and we have 
	\begin{align}
		\IPMP(d(z);N) \coloneqq \left\{ \begin{matrix} \RFMP(d(z);N), &\text{if LRFMP is used,}\\ \ROFMP_S(d(z);N), &\text{if LROFMP is used.} 	\end{matrix} \right.
	\end{align}
For technical reasons, we use $\ROFMP_S(d(z);N)$ which is the product of $\ROFMP(d(z);N)$ from \cref{def:ROFMP(d;N)} and a spline to avoid neighbourhoods of previously chosen basis functions from the respective trial function class. Let $\varepsilon$ denote the size of such a neighbourhood. To avoid a neighbourhood of a chosen $z$, such a spline is given by
\begin{align}
		\left( S_{z} \right)|_{[0,\varepsilon]} &\equiv 0,\qquad 
		\left(S_{z}\right)|_{(\varepsilon,2\varepsilon)}(\tau) &\coloneqq 10\left(\frac{\tau}{\varepsilon} -1\right)^3 -15\left(\frac{\tau}{\varepsilon}-1\right)^4 + 6\left(\frac{\tau}{\varepsilon}-1\right)^5
		\qquad 	\textrm{and} \qquad
	\left( S_{z} \right)|_{[2\varepsilon,2]} &\equiv 1.
	\end{align}
In the $N$-th iteration, this yields  
	\begin{align}
		\ROFMP_S(d(z);N) = \ROFMP(d(z);N) \cdot \prod_{n=1}^{N} S_{z^{(n)}}\left( \left\| z - z^{(n)} \right\|^2\right),
		\label{eq:objfun:spline}
	\end{align}
where $\tau = \| z - z^{(n)} \|^2$ denotes the distance between the current value $z$ and previously chosen values $z^{(n)}$. Note that we only have to consider the spline if $d(z)$ and $d(z^{(n)})$ are from the same trial function class. Then, for each truly infinite trial function class, we solve the maximization problem $\IPMP(d(z);N) \to \max!$ in each iteration $N$.\\
The replacement character $d(z)$ stands for either $g^{(k,L)}(R,\cdot)$ with $z=R(c,\alpha,\beta,\gamma) \in \ball^4$ or for $K(x,\cdot)$ and $W(x,\cdot)$ with $z=x\in \ball$. Hence, we maximize $\IPMP(d(z);N)$ with respect to the characteristic parameter vector $z$ of the trial functions SL, APK and APW. Thus, we have to model the corresponding constraints as well. For the SLs, we have in practice
	\begin{align}
		(c,\alpha,\beta,\gamma) \in [-1,1] \times [0,2\pi] \times [0,\pi] \times [0,2\pi].
	\end{align}
For the APKs and APWs, we obtain
	\begin{align}
		|z| < 1.
	\end{align}
In practice, we solve these optimization problems using the NLOpt library, see \cite{NLOpt2019}. In particular, as it is advised there, we solve them in a 2-step-optimization procedure. That means, we first determine a global solution (with a derivative-free method) and, then, refine this using a gradient-based local method. We prefer to choose the local solution as our candidate of the respective trial function class. Note that, in the case of problems with a solver, we can also use the global solution. If the global solver needs a sensible starting solution, we can include a selection of SLs, APKs or APWs in the starting dictionary as well. However, these should not have a major impact on the learnt dictionary.\\
In \cite{Micheletal2018-1}, we proposed the use of certain additional features to guide the learning process. Though the features proved to be helpful in certain learning settings, from our experience, using a 2-step optimization, i.\,e. solving the described optimization problems first globally and then locally, as well as using more diverse trial function classes remedies the urgent need of some rather manual features. Nonetheless, some of them like an iterative application of the learnt dictionary (i.\,e. allowing only the first $N$ dictionary elements in the $N$-th iteration of the IPMP when the learnt dictionary is used) is in particular helpful when we have to balance the tradeoff between numerical accuracy and runtime. \\
As an overview, we provide a hierarchy of two pseudocodes for the LIPMP algorithms. \cref{alg:PseudocodeLIPMP} describes the determination of the dictionary of candidates and is integrated into \cref{alg:PseudocodeIPMP} which is the main routine. Note that $[\cdot^{\mathrm{s}}]$ stands for the respective trial function class from the starting dictionary. Further, note that we did not include the restart procedure of the iterated (L)ROFMP algorithm in order to concentrate on the novel aspects of the add-on. The restart technique was described in \cite{Micheletal2016-1,Telschow2014,Schneider2020}.

	\begin{algorithm}[htbp]
		\KwData{$R^N \in \real^\ell,\ N \in \nat_0,\ \dic^\mathrm{s} $}
		\KwResult{choice of next basis element $d_{N+1}$}
		\begin{itemize}
			\item[(1)] determine \\
				$$d_{N+1}^\SH = \mathrm{argmax}_{d\in [\cdot^{\mathrm{s}}]_\SH}  \IPMP(d;N);$$
			\item[(2)] compute starting choices $d_{N+1}^{(\mathrm{s},\bullet)}(z),\ \bullet \in \{\APK,\ \APW,\ \SL \}$:
				$$d_{N+1}^{(\mathrm{s},\bullet)} = \mathrm{argmax}_{d(z)\in [\cdot^{\mathrm{s}}]_\bullet}  \IPMP(d(z);N);$$
			\item[(3)] compute global choices $d_{N+1}^{(\mathrm{g},\bullet)}(z),\ \bullet \in \{\APK,\ \APW,\ \SL \}$ using starting choices as starting points in the optimization processes:
					$$d_{N+1}^{(\mathrm{g},\APK)} = \mathrm{argmax}_{d(z)\in \left[\ball\right]_\APK}  \IPMP(d(z);N);$$
					$$d_{N+1}^{(\mathrm{g},\APW)} = \mathrm{argmax}_{d(z)\in \left[\ball\right]_\APW}  \IPMP(d(z);N);$$
					$$d_{N+1}^{(\mathrm{g},\SL)} = \mathrm{argmax}_{d(z)\in \left[\ball^4\right]_\SL}  \IPMP(d(z);N);$$
			\item[(4)] compute local solutions $d_{N+1}^{(\mathrm{l},\bullet)}(z),\ \bullet \in \{\APK,\ \APW,\ \SL \}$ using global choices as starting points in the optimization processes:
				$$d_{N+1}^{(\mathrm{l},\APK)} = \mathrm{argmax}_{d(z)\in \left[\ball\right]_\APK}  \IPMP(d(z);N);$$
				$$d_{N+1}^{(\mathrm{l},\APW)} = \mathrm{argmax}_{d(z)\in \left[\ball\right]_\APW}  \IPMP(d(z);N);$$
				$$d_{N+1}^{(\mathrm{l},\SL)} = \mathrm{argmax}_{d(z)\in \left[\ball^4\right]_\SL}  \IPMP(d(z);N);$$
			\item[(5)] choose 
				$$d_{N+1}^{\phantom{(\mathrm{l},\bullet)}} = \mathrm{argmax}_{\hat{d} \in \Can} \IPMP\left(\hat{d};N\right),$$
				$$\Can \coloneqq \left\{\left. d_{N+1}^\SH,\ d_{N+1}^{(\mathrm{s},\bullet)},\ d_{N+1}^{(\mathrm{g},\bullet)},\,\ d_{N+1}^{(\mathrm{l},\bullet)}\ \right|\ \right. $$
				$$\left. \qquad \qquad \qquad \qquad \vphantom{d_{N+1}^\SH} \bullet \in \{\APK,\ \APW,\ \SL \}  \right\};$$
		\end{itemize}
	return $d_{N+1}$\;
	\vspace*{\baselineskip}
	\caption{Pseudo-code for choosing a dictionary element from $\dicinf$ for the LIPMP algorithms.}
	\label{alg:PseudocodeLIPMP}
	\end{algorithm}
	\begin{algorithm}[htbp]
		\KwData{$y \in \real^\ell$}
		\KwResult{approximation $f_N$}
			initialization: $\dic^\mathrm{s},\ f_0 \equiv 0,\ R^0 = y - T_\daleth f_0 = y$ \;
			$N = 0$\;
		\While{(stopping criteria not fulfilled)}{
			determine $d_{N+1}$ via \cref{alg:PseudocodeLIPMP} with $\IPMP(\cdot;N)$\;
			compute
			\newline
				\begin{tabular}{rl}			
					$\alpha_{N+1}$ & via \cref{def:alphan} or \cref{def:alphanN} \\ &\\
					$R^{N+1}$ & via \cref{def:RN} or \cref{def:RNO} \\ &\\
					$N$ & be increased by 1;\\&\\
					and & process the newly chosen basis element 
			\end{tabular}
		}
			return approximation via \cref{eq:fn} or \cref{eq:fnN}\;
	\vspace*{\baselineskip}
	\caption{Pseudo-code for the LIPMP algorithm.}
	\label{alg:PseudocodeIPMP}
	\end{algorithm}

\subsection{Summary of properties}
Before we present our numerical results, we summarize some properties of the LIPMP algorithms. For more details on these aspects, we refer the reader to \cite{Schneider2020}.\\
The LIPMP algorithms can be used as standalone approximation algorithms or as auxiliary algorithms to determine a finite dictionary automatically. Hence, they yield an approximation of the inverse problem as well as a learnt dictionary for this problem. \\
Relating them to other dictionary learning algorithms, see e.\,g. \cite{Aharonetal2006,Brucksteinetal2009,Enganetal1999,Enganetal1999-1,Pruente2008,Rubinsteinetal2010}, we find that they similarly solve a doubled minimization problem. Furthermore, if the criteria for dictionary learning presented in this literature are re-interpreted for the case of inverse problems and continuous approximations, the LIPMP algorithms also fulfil such aspects. \\
By construction, they inherit the convergence results of the IPMP algorithms (see the literature mentioned above). In particular, we have convergence of the approximation for infinitely many iterations in the LRFMP algorithm to the solution of the regularized normal equation. \\
Further, a sensible criterion for a sequence of well-working dictionaries is that, in the limit, the solution of the regularized normal equation can be represented. The iteration process of the LIPMP algorithms naturally defines a nested sequence of learnt dictionaries (one additional sequence element per iteration). Combining this with the convergence results, we have that, in the case of the LRFMP algorithm, the sequence of learnt dictionaries is a sequence of well-working dictionaries. 

\section{Numerical results}
\label{sect:numresults}
We first summarize the general setting of the experiments. Then we show results for comparing a manually chosen and a learnt dictionary in the IPMP algorithms. At last, we show the results of the LIPMP algorithms as standalone approximation methods. The results are presented in more detail in \cite{Schneider2020}. Note that our test scenarios here shall serve as proofs of concept in the sense that the main features of the add-on are demonstrated. In a continuing project, we will investigate the behaviour for more realistic data and for other applications.

\subsection{Experiment setting}
For a longer description of the experiment setting see \cite{Schneider2020}. As data, we use the EGM2008, see e.\,g. \cite{Pavlisetal2012}, as well as GRACE data from May 2008 evaluated on a regularly distributed Reuter grid of 12684 points. For the GRACE data, we utilize the degree 3 to 60 from the arithmetic mean of the Level 2 Release 05 provided by the GFZ, JPL and UTCSR. Further, we smooth the data with a Cubic Polynomial Spline of order 5, see \cite{Schreiner1996,Freedenetal1998,Fengleretal2006} to remove the North-South striping.\\
For the relative RMSE, we utilize this data on an equi-angular Driscoll-Healy grid \cite{DriscollHealy1994,Michel2013} with 65341 points (confer \cref{fig:dataandgrid}). The data is modelled on a 500 km satellite height and is perturbed with $5\%$ Gaussian noise, such that we have perturbed data $y^\delta$ given by 
	\begin{align}
		y^\delta_i = y_i \cdot \left( 1 + 0.05 \cdot \varepsilon_i\right)
	\end{align}
for the unperturbed data $y_i$ and a Gaussian distributed random number $\varepsilon_i$. \\
The algorithms terminate if the relative data error falls below the noise level or if 1000 iterations are reached. We implemented the iterated (L)ROFMP algorithm and restarted the procedure after 100 iterations. We choose the tested regularization parameter that minimized the relative RMSE if the relative data error reached the noise level at termination.\\
The optimization problems are solved by the ORIG\_DIRECT\_L (globally) and the SLSQP (locally) algorithms from the NLOpt. As it is advised, we narrow the constraints by $10^{-8}$. Further, we set some termination criteria of the optimization procedures: we limit the absolute tolerance of the change of the objective function between two successive iterates as well as the tolerance between the iterates themselves to $10^{-8}$. Moreover, we allow 5000 function evaluations and 200 seconds for each optimization.\\
With respect to the SLs, APKs and APWs, we forbid to choose two trial functions of the same type which are as close as $5 \cdot 10^{-4}$ or closer in one (L)ROFMP step. Further, we use the same regularization parameter for learning and applying the learnt dictionary. The regularization parameter is constant unless anything different is stated. We apply the dictionary learning iteratively (confer \cite{Micheletal2018-1}).\\
As the starting dictionary, we use 
	\begin{align}
		\left[N^\mathrm{s}\right]_\SH &= \left\{ Y_{n,j}\ \left|\ n=0,...,100;\ j=-n,...,n \right.\right\}\\
		\left[S^\mathrm{s}\right]_\SL &= \left\{ \left. g^{(k,5)}\left(\left(c,A(\alpha,\beta,\gamma)\ethree\right),\cdot \right) \right| \right. \\
		& \qquad \quad c\in\left\{\frac{\pi}{4},\frac{\pi}{2}\right\},\ \alpha \in \left\{0,\frac{\pi}{2},\pi,\frac{3\pi}{2}\right\}, \beta\in\left\{0,\frac{\pi}{2},\pi\right\},\ \gamma\in\left\{0,\frac{\pi}{2},\pi,\frac{3\pi}{2}\right\}, \left. \vphantom{g^{(k,5)}} k=1,...,36\right\}\\
		\left[B_K^\mathrm{s}\right]_\APK &= \left\{ \frac{K(x,\cdot)}{\|K(x,\cdot)\|_{\Lp{2}}}\ \left|\ |x| = 0.94,\ \frac{x}{|x|} \in X^\mathrm{s} \right. \right\}	\\
		\left[B_W^\mathrm{s}\right]_\APW &= \left\{ \frac{W(x,\cdot)}{\|W(x,\cdot)\|_{\Lp{2}}}\ \left|\ |x| = 0.94,\ \frac{x}{|x|} \in X^\mathrm{s} \right.\right\}\\[.5\baselineskip]
		\dic^{\mathrm{s}} &= \left[N^\mathrm{s}\right]_{\SH} \cup \left[S^\mathrm{s}\right]_{\SL} \cup  \left[B_K^\mathrm{s}\right]_{\APK} \cup \left[B_W^\mathrm{s}\right]_{\APW}
	\end{align}
with a regularly distributed Reuter grid $X^\mathrm{s}$ of $123$ grid points. Thus, the starting dictionary contains 13903 trial function.\\
The experiments that include the LIPMP algorithms ran on a node of 48 GB RAM with 12 CPUs.

\subsection{Learning a dictionary}
\label{ssect:Comp}
\begin{table*}
\begin{minipage}{\textwidth}
\begin{center}
\footnotesize
\begin{tabular}{lllllllll}
\hline
data 			& algorithm/					&regularization 								& size of 			&completed		& maximal 			& relative  		& relative 			& CPU-runtime\\ 
					& setting						&parameter										& dictionary		&iterations		& degree				& RMSE			& data error		& in h	\\
\hline
EGM2008 	& $\RFMP^*$ 				&	$10^{-9}\|y\|_{\real^\ell}$		& 95152				& 957					&		25					& 0.000466	& 0.049998		& 514.03\\
					&$\RFMP^{**}$ 			&	$10^{-9}\|y\|_{\real^\ell}$		& $\leq 637$		& 662					&		46					& 0.000471 	& 0.049999 		& 507.22\\ 
					&$\RFMP^{***}$ 		&	$10^{-6}\|y\|_{\real^\ell}/N$	& $\leq 670$		& 670					&		34					& 0.000447 	& 0.050000 		& 533.87\\ 
					&$\RFMP^{****}$ 		&	$10^{-9}\|y\|_{\real^\ell}$		& $\leq 684$		& 734					&		41					& 0.000484 	& 0.049999 		& 129.57\\ 
\hline
EGM2008	&$\ROFMP^*$ 			&	$10^{-9}\|y\|_{\real^\ell}$		& 95152				& 766					&		25					& 0.000463	& 0.049999		& 561.76\\
					&$\ROFMP^{**}$		&	$10^{-9}\|y\|_{\real^\ell}$		&	$\leq 550$		& 577					&		38					& 0.000467	& 0.049998 		& 665.76\\ 
					&$\ROFMP^{***}$		&	$10^{-6}\|y\|_{\real^\ell}/N$	&	$\leq 686$		& 701					&		35					& 0.000452	& 0.049999 		& 840.06\\ 
					&$\ROFMP^{****}$	&	$10^{-9}\|y\|_{\real^\ell}$		&	$\leq 600$		& 621					&		46					& 0.000468	& 0.049998 		& 386.08\\ 
\hline	
GRACE		&$\RFMP^*$ 				&	$10^{-4}\|y\|_{\real^\ell}$		& 95152				& 393					&		25					& 0.000340	& 0.049997		& 522.09\\
					&$\RFMP^{**}$			&	$10^{-4}\|y\|_{\real^\ell}$		&	$\leq 384$		& 483					&		32					& 0.000335	& 0.049999 		& 341.16\\ 
					&$\RFMP^{***}$			&	$10^{-1}\|y\|_{\real^\ell}/N$	&	$\leq 352$		& 349					&		28					& 0.000311	& 0.049994 		& 284.85\\ 
					&$\RFMP^{****}$		&	$10^{-4}\|y\|_{\real^\ell}$		&	$\leq 479$		& 535					&		39					& 0.000344	& 0.049994 		& \phantom{0}90.53\\ 
\hline
GRACE		&$\ROFMP^*$ 			&	$10^{-4}\|y\|_{\real^\ell}$		& 95152				& 274					&		25					& 0.000328	& 0.049989		& 528.67\\
					&$\ROFMP^{**}$		&	$10^{-4}\|y\|_{\real^\ell}$		&	$\leq 303$		& 306					&		34					& 0.000330 	& 0.049998 		& 372.31\\ 
					&$\ROFMP^{***}$		&	$10^{-3}\|y\|_{\real^\ell}/N$	&	$\leq 292$		& 290					&		33					& 0.000374	& 0.049998 		& 358.64\\ 
					&$\ROFMP^{****}$	&	$10^{-4}\|y\|_{\real^\ell}$		&	$\leq 278$		& 278					&		26					& 0.000322 	& 0.049996 		& 182.19\\ 
\hline
\end{tabular}
\normalfont
\caption{Comparison of a manually chosen and diverse learnt dictionary. Confer \cref{ssect:Comp}. The IPMP$^{*}$ algorithm uses the manually chosen dictionary, the IPMP$^{**}$ the learnt dictionary, the IPMP$^{***}$ the non-stationary learnt dictionary and the IPMP$^{****}$ the learnt-without-Slepian-functions dictionary. All learnt dictionaries are iteratively applied. The maximal degree is the maximal SH degree included in the used dictionary. $N \in \nat$ stands for the iterations.}
\label{tab:T12} 
\end{center}
\end{minipage}
\end{table*}
We compare the learnt dictionary with a manually chosen dictionary which is similar to those in previous publications, see e.\,g. \cite{Telschow2014}: 
	\begin{align}
		\left[N^\mathrm{m}\right]_\SH &= \left\{\left. Y_{n,j}\ \right|\ n=0,...,25;\ j=-n,...,n \right\}\\
		\left[S^\mathrm{m}\right]_\SL &= \left\{ \left. g^{(k,5)}\left(\left(c,A(\alpha,\beta,\gamma)\ethree\right),\cdot \right) \right| \right. \\
		& \qquad \quad c\in\left\{\frac{\pi}{4},\frac{\pi}{2}\right\},\ \alpha \in \left\{0,\frac{\pi}{2},\pi,\frac{3\pi}{2}\right\}, \beta\in\left\{0,\frac{\pi}{2},\pi\right\},\ \gamma\in\left\{0,\frac{\pi}{2},\pi,\frac{3\pi}{2}\right\},\left. \vphantom{g^{(k,5)}} k=1,...,36\right\}
	\end{align}
	\begin{align}
		\left[B_K^\mathrm{m}\right]_\APK &= \left\{ \frac{K(x,\cdot)}{\|K(x,\cdot)\|_{\Lp{2}}}\ \left|\ |x| \in Z,\ \frac{x}{|x|} \in X^\mathrm{m} \right. \right\}		\\
		\left[B_W^\mathrm{m}\right]_\APW &= \left\{ \frac{W(x,\cdot)}{\|W(x,\cdot)\|_{\Lp{2}}}\ \left|\ |x| \in Z,\ \frac{x}{|x|} \in X^\mathrm{m} \right.\right\}\\[.5\baselineskip]
		\dic^{\mathrm{m}} &= \left[N^\mathrm{m}\right]_{\SH} \cup \left[S^\mathrm{m}\right]_{\SL} \cup  \left[B_K^\mathrm{m}\right]_{\APK} \cup \left[B_W^\mathrm{m}\right]_{\APW}
	\end{align}
with a regularly distributed Reuter grid $X^\mathrm{m}$ of $4551$ grid points and 
	\begin{align}
		Z= \{0.75,\ 0.80,\ 0.85,\ 0.89,\ 0.91,  0.93,\ 0.94,\ 0.95,\ 0.96,\ 0.97\}.
	\end{align}
All in all, the manually chosen dictionary contains $95152$ trial functions.
We undertake this comparison because it is most sensible as we have explained in \cite{Micheletal2018-1}: a comparison with the best dictionary of a sensibly large set of random dictionaries cannot seriously be put into practice due to high memory demand and a long runtime. Note that, in some of the literature on the IPMP algorithms mentioned before, the approaches have been compared to traditional methods like splines in previous publications as well. Further note that, due to the size of the manually chosen dictionary, the respective tests ran on a node with 512 GB RAM and 32 CPUs.\\
In \cref{tab:T12}, we see a summary of the results of the experiments. We compare the IPMP algortihms with the manually chosen dictionary (*), the learnt dictionary (**), a learnt dictionary when using the non-stationary regularization parameter $\lambda_N = \lambda_0 \cdot \|y\|_\real^\ell/N$ for the iteration $N\in\nat$ (``non-stationary learnt"; ***), and a learnt dictionary where only the SHs, APKs and APWs were considered (``learnt-without-Slepian-functions"; ****). We give the regularization parameter, the size of the dictionary, the number of completed iterations, the maximal SH degree included in the dictionary, the relative data error and RMSE at termination and the needed CPU-runtime in hours. Note that the size of the learnt dictionaries is given as a ``less or equal than" value since elements may be contained multiple times.\\
We see the following aspects:
	\begin{compactitem}
		\item Due to our termination criteria, the relative data error was at the noise level in all cases.
		\item In all experiments, the relative RMSE is about the same size. In comparison to \cite{Micheletal2018-1}, we conclude that the IPMP algorithms produce better results when more trial function classes are used. Further, the learnt dictionary yields similar results. In Figures \ref{fig:supp2} and \ref{fig:supp3}, we also see that, in all cases, the remaining errors lie within regions of higher local structures, i.\,e.\ the Andean region, the Himalayas and the Pacific Ring of Fire in the case of EGM2008 data as well as the Amazon basin in the case of the GRACE data. These detail structures cannot be represented beause of the noise and the ill-posedness.
		\item The non-stationary learnt as well as the learnt-without-Slepian-functions dictionary produce similar results with respect to the relative RMSE such that these settings could be explored in future research. However, to evaluate the influence of the non-stationary regularization parameter on the approximation, the number of tests is here to low.
		\item The learnt dictionary is less than $1\%$ of the size of the manually chosen dictionary.
		\item The maximal SH degree of the learnt dictionaries is a truly learnt degree.
		\item For the LRFMP algorithm, the CPU-time needed for learning and applying the learnt dictionary is lower or similar than applying the manually chosen dictionary. In particular, without the Slepian functions, the needed CPU-time is much smaller.
		\item For the LROFMP algorithm, there exist settings which have a smaller runtime as well. In particular, for the GRACE data, this is always the case. Similarly, the learnt-without-Slepian-functions dictionary is also learnt in a much shorter time for the EGM2008 data. However, there are two cases for the EGM2008 data where the runtime is higher than for the manually chosen dictionary. This could be caused by the non-stationary regularization parameter, the orthogonality procedure itself and / or the use of the Slepian functions.
		\item At last, the learning of a dictionary and its use in the respective IPMP has a much lower storage demand than in the case of the manually chosen dictionary.
	\end{compactitem}
Note that, in \cite{Schneider2020}, we also presented a first approach to learn a dictionary and apply it to unseen test data. The results there were satisfactory but should be improved in the future.

\subsection{The LIPMP algorithms as standalone approximation methods}
\label{ssect:STD}
\begin{table*}
\begin{minipage}{\textwidth}
\begin{center}
\footnotesize
\begin{tabular}{lllllll}
\hline
experiment 							& algorithm		&regularization 							&completed		& maximally	learnt   	& relative  			& relative 			\\ 
												&							&parameter									&iterations		& SH degree					& RMSE				& data error			\\
\hline
approximation						& LRFMP			& $\phantom{5\cdot \ }\ 10^{-9}\|y\|_{\real^\ell}$ 	& 1000				& 75 								& 0.000249 	& 0.075293\\
												& LROFMP			& $\phantom{5\cdot \ }\ 10^{-9}\|y\|_{\real^\ell}$ 	& 1000				& 83 								& 0.000253 	& 0.075210\\
\hline
downward continuation,	& LRFMP			& $\phantom{5\cdot \ }\ 10^{-9}\|y\|_{\real^\ell}$ 	& 637					& 46								& 0.000471 	& 0.049995\\
regular grid							& LROFMP			& $\phantom{5\cdot \ }\ 10^{-9}\|y\|_{\real^\ell}$ 	&	550 				& 38								& 0.000465 	& 0.049998\\
\hline
downward continuation,	& LRFMP			& $5\cdot 10^{-9}\|y\|_{\real^\ell}$ 	& 975  				& 50 								& 0.000472	& 0.050000\\
irregular grid							& LROFMP			& $5\cdot 10^{-9}\|y\|_{\real^\ell}$ 	& 983 				& 50 								& 0.000521	& 0.049999\\
\hline			
contrived data						& LROFMP			& $\phantom{5\cdot \ }\ 10^{-8}\|y\|_{\real^\ell}$	& 24					& 9 									& 0.0000076 & 0.049975\\
\hline
\end{tabular}
\normalfont
\caption{The LIPMP algorithms as standalone approximation methods. Confer \cref{ssect:STD}. }
\label{tab:STD}
\end{center}
\end{minipage}
\end{table*}
\begin{figure}[tb]
	\begin{subfigure}{\textwidth}	
			\centering
			\includegraphics[width=.45\textwidth]{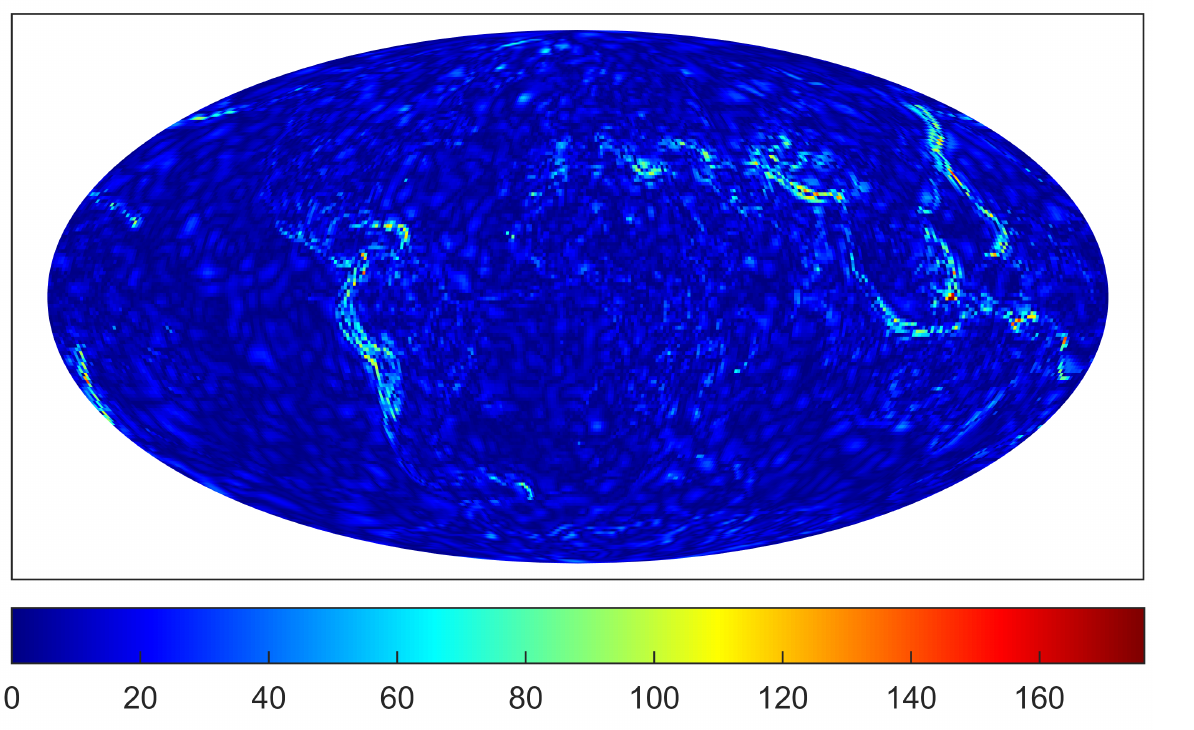}
			\includegraphics[width=.45\textwidth]{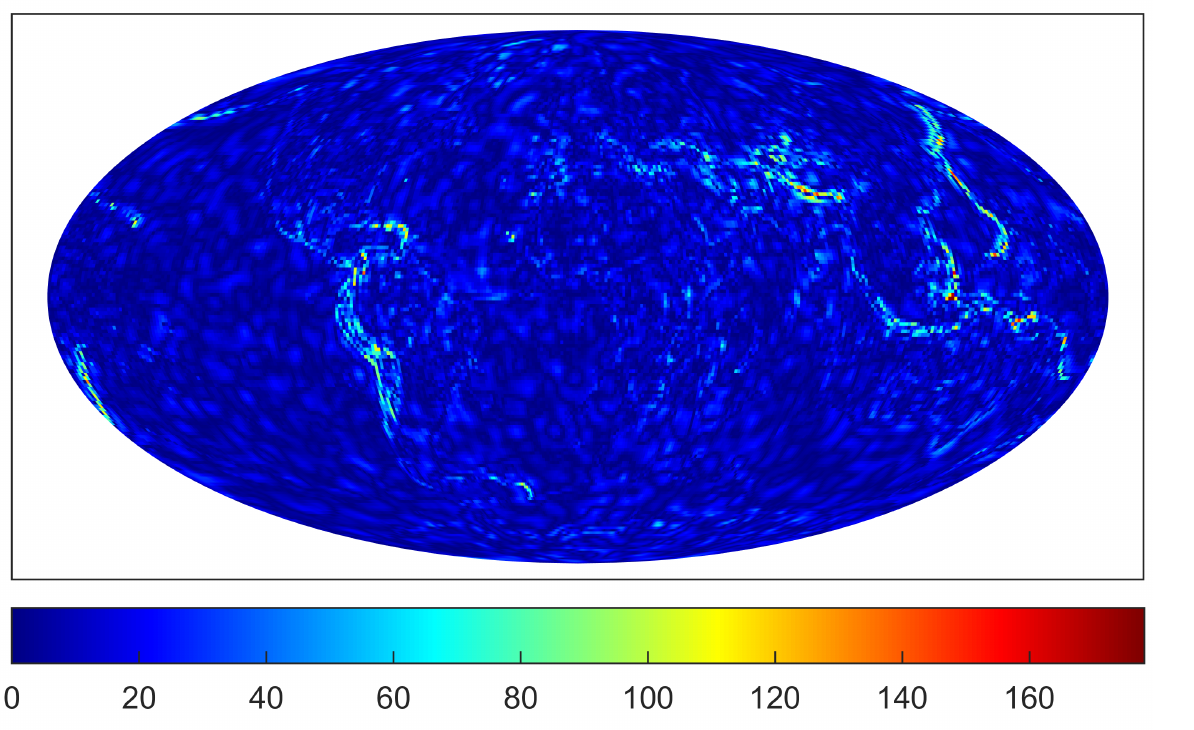}
			\caption{Absolute approximation error obtained by the LRFMP (left) and the LROFMP (right) algorithm for the approximation of surface data. All values in $\mathrm{m}^2/\mathrm{s}^2$.}
			\label{fig:T4}
	\end{subfigure}
	\begin{subfigure}{\textwidth}	
			\centering
			\includegraphics[width=.45\textwidth]{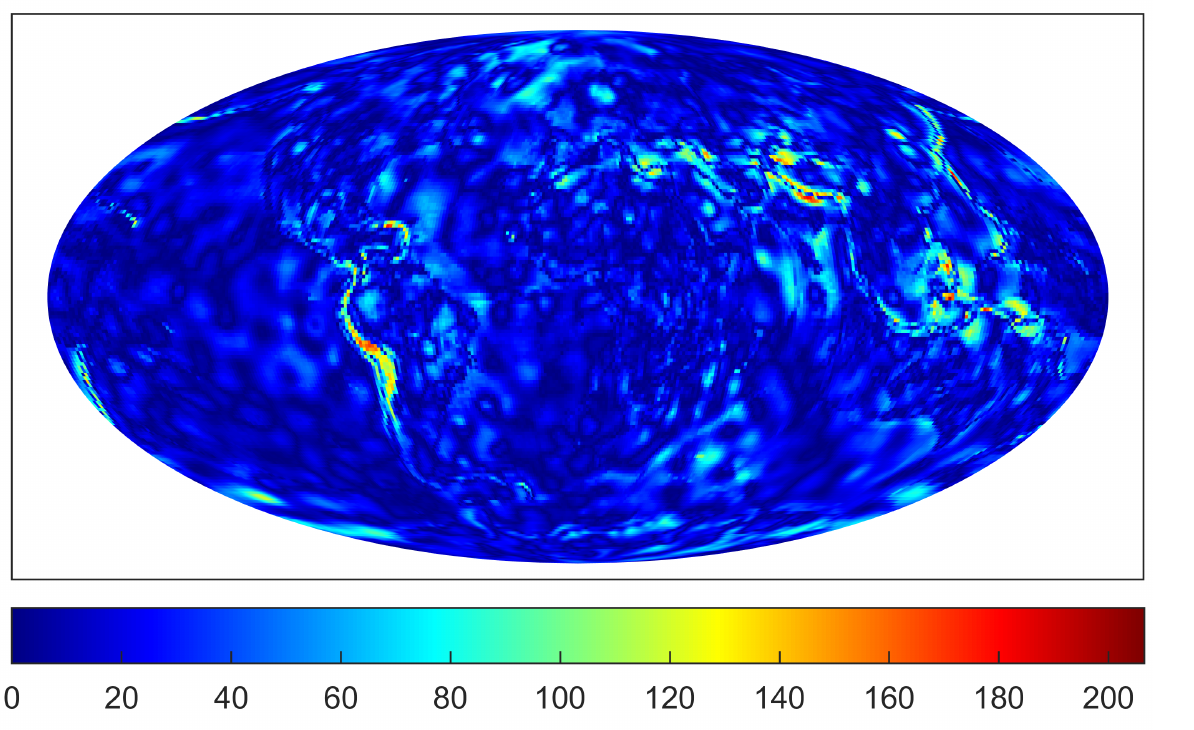}
			\includegraphics[width=.45\textwidth]{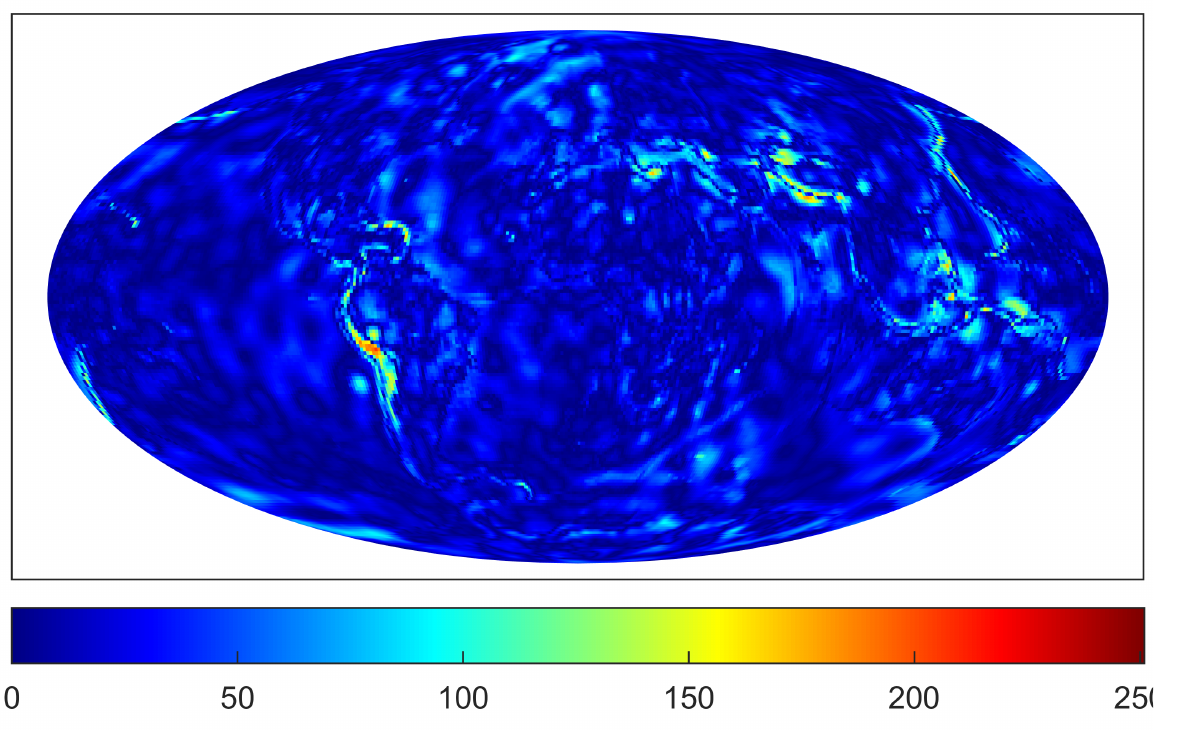}
			\caption{Absolute approximation error obtained by the LRFMP (left) and LROFMP (right) algorithm for the downward continuation of satellite data form a regularly distributed grid. The scales are adapted for a better comparison with \cref{fig:supp2}. All values in $\mathrm{m}^2/\mathrm{s}^2$.}
			\label{fig:T12}
	\end{subfigure}
\caption{Results of the LIPMP algorithms as standalone approximation algorithms. Confer \cref{ssect:STD}.}
\end{figure}
\begin{figure}[tb]
	\begin{subfigure}{\textwidth}
			\centering
			\includegraphics[width=.45\textwidth]{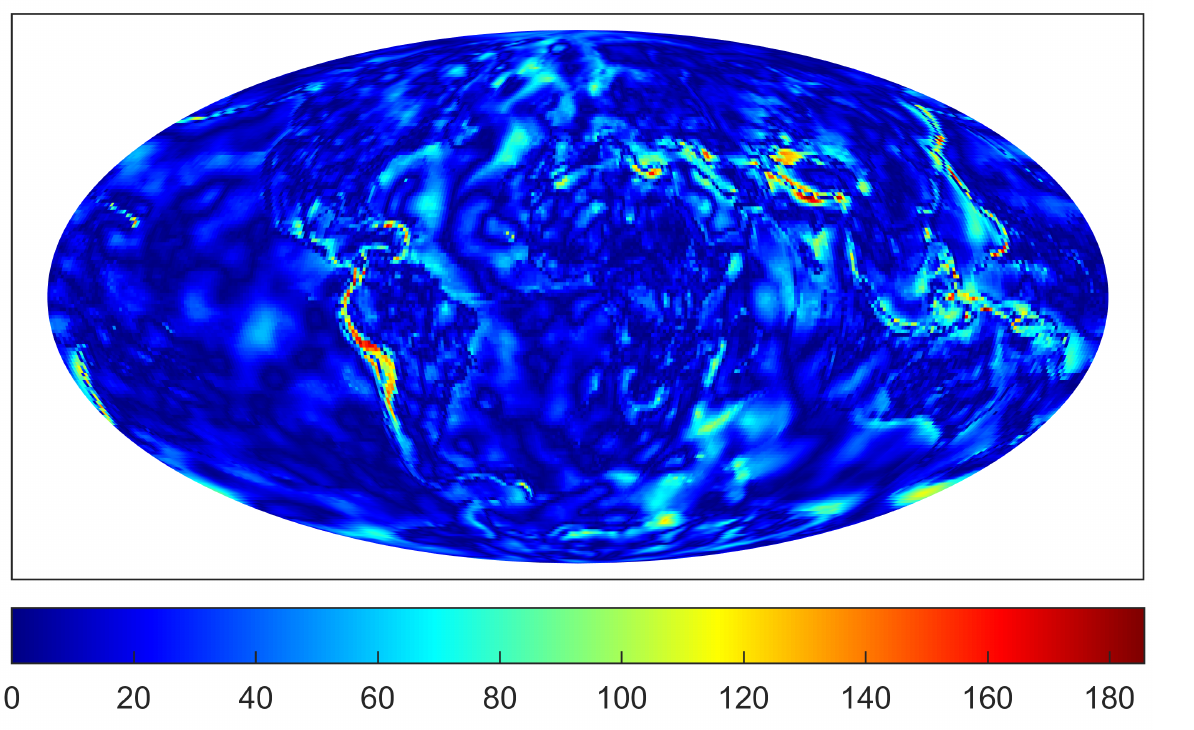}
			\includegraphics[width=.45\textwidth]{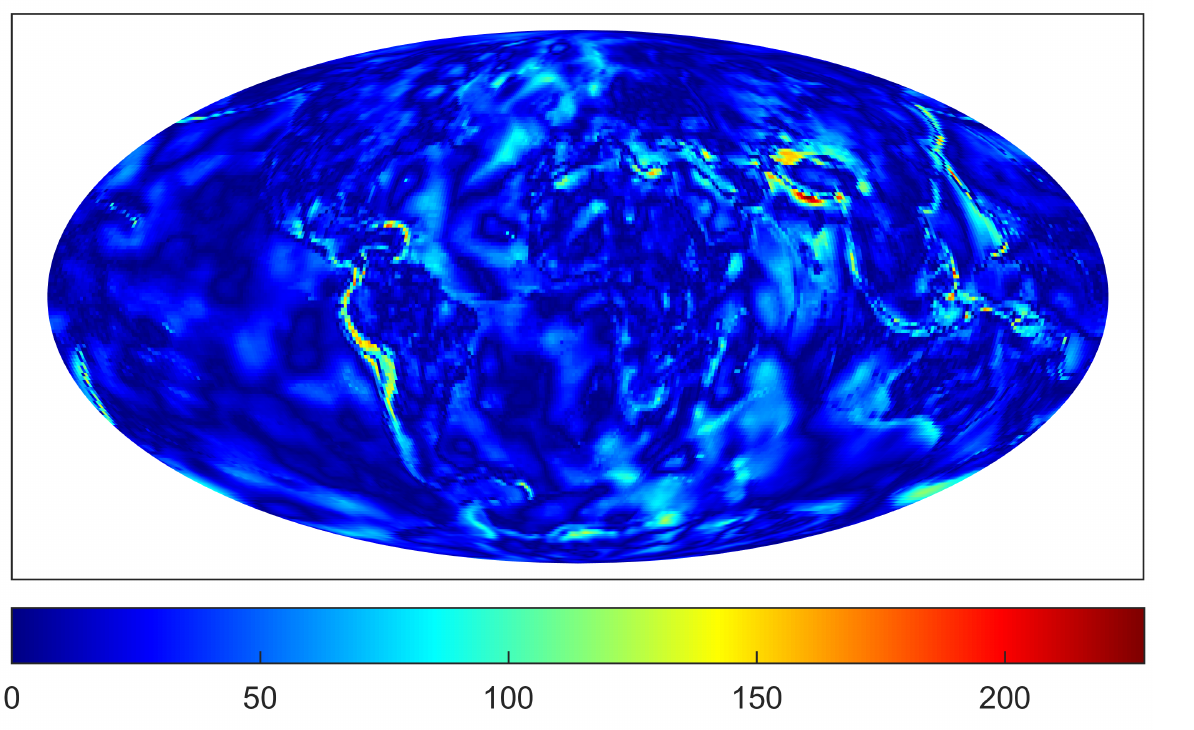}
			\caption{Absolute approximation error obtained from the LRFMP (left) and the LROFMP (right) algorithm for the downward continuation of satellite data from an irregularly distributed grid. All values in $\mathrm{m}^2/\mathrm{s}^2$.}
			\label{fig:T7}
	\end{subfigure}
	\begin{subfigure}{\textwidth}
		\centering
		\includegraphics[width=.75\textwidth]{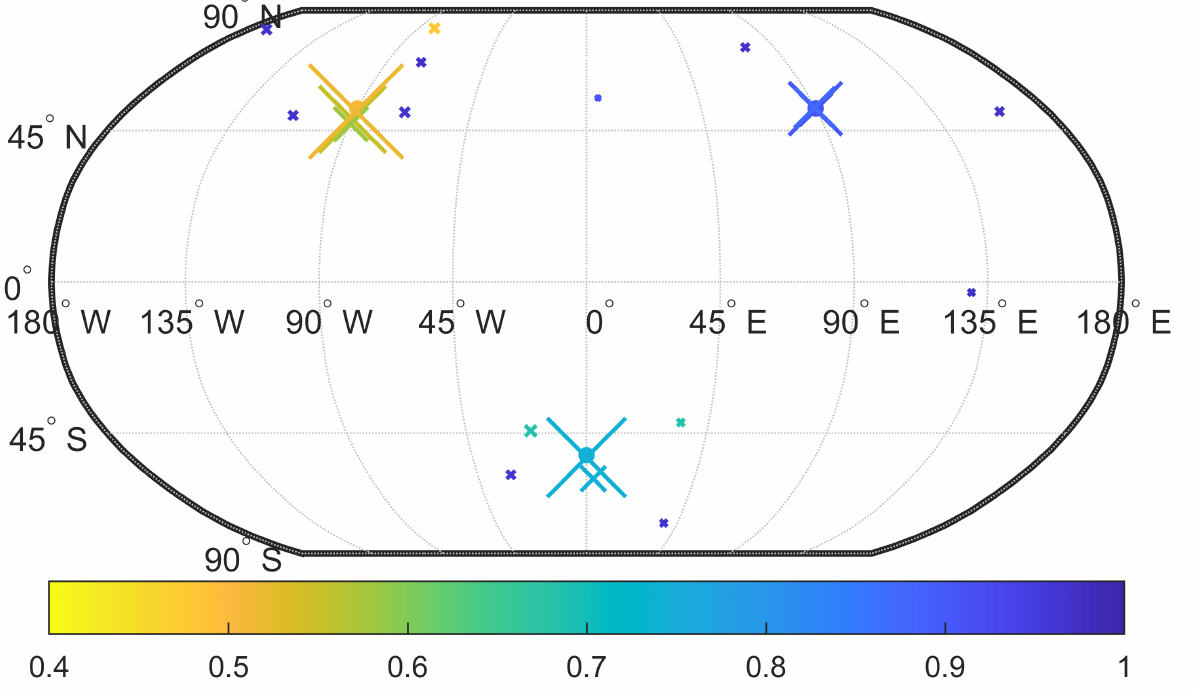}
		\caption{Given and chosen APKs for perturbed contrived data. The dots stand for the given solution, the crosses give the chosen APKs. The colour represents the respective scale and the sizes visualize the absolute value of the respective chosen coefficients.}
		\label{fig:T5-5}
	\end{subfigure}
\caption{Further results of the LIPMP algorithms as standalone approximation algorithms. Confer \cref{ssect:STD}.}
\end{figure}
So far, we used the LIPMP algorithms to learn a dictionary. This dictionary was used afterwards in the corresponding IPMP algorithm to construct a regularized approximation of the gravitational potential. However, the learning algorithm itself incorporates the maximization of the same objective function which also occurs in the IPMP. Indeed, we observed that the LIPMP already produces a very good approximation without an additional run of the IPMP. This shall be demonstrated in this section. For this purpose, we consider the approximation of surface data, the downward continuation of regularly and irregularly distributed satellite data by both LIPMP algorithms. Moreover, we verify the downward continuation of contrived data by the LROFMP algorithm. The irregularly distributed grid has already been used in \cite{Micheletal2014} and simulates a denser data distribution on the continents. It is given in \cref{fig:dataandgrid} and includes 6968 grid points. The contrived data consists of 3 SHs and APKs, respectively: 
\begin{align}
		f 	&= Y_{9,5}\left(\cdot\right) + Y_{5,5}\left(\cdot\right) + Y_{2,0}\left(\cdot\right) + \widetilde{K}\left(x\left(0.5,\frac{3\pi}{2},\frac{\pi}{4}\right),\cdot\right) + \widetilde{K}\left(x\left(0.75,2\pi,-\frac{\pi}{4}\right),\cdot\right) + \widetilde{K}\left(x\left(0.9,\frac{\pi}{2},\frac{\pi}{4}\right),\cdot\right) 
	\end{align}
where the notation $x(r,\lon,\cos(\lat))$ with the radius $r$, the longitude $\lon$ and the latitude $\lat$ is used. $\widetilde{K}$ stands for the $\Lp{2}$-nor\-mal\-ized APKs. The data is again perturbed by $5\%$ Gaussian noise. Correspondingly, the starting dictionary (for the test with contrived data only) is given as
\begin{align}
		\left[N^\mathrm{s}\right]_\SH &= \left\{\left. Y_{n,j}\ \right|\ n=0,...,10;\ j=-n,...,n  \right\}\\
		\left[B_K^\mathrm{s}\right]_\APK &= \left\{ \frac{K(x,\cdot)}{\| K(x,\cdot)\|_{\Lp{2}}}\ \left|\ |x| = 0.94,\ \frac{x}{|x|} \in X^\mathrm{s} \right. \right\}\\
		\dic^{\mathrm{s}} &= \left[N^\mathrm{s}\right]_{\SH} \cup\left[B_K^\mathrm{s}\right]_{\APK},
	\end{align}
where $X^\mathrm{s}$ is a regularly distributed Reuter grid of $6$ points. Then the APKs included in the contrived data are not contained in the starting dictionary. We allow a maximum of 100 iterations here because the data consists of only six trial functions. Due to the orthogonality procedure, we assume that the LROFMP algorithm is more suited to obviously distinguish the SHs and APKs.\\
In \cref{tab:STD}, we give an overview of the results. The type of experiment is abbreviated: ``approximation" stands for the experiment where no satellite height is included, ``downward continuation, (ir-)\-regular grid" stands for the use of 500 km satellite height and an (ir-)\-regularly distributed grid and ``contrived data" is self-explanatory. Further, we state the regularization parameter, the completed iterations, the respective maximally learnt SH degrees and the relative RMSEs as well as data errors. In \crefrange{fig:T4}{fig:T7}, we see the absolute approximation errors obtained in the different experiments. \cref{fig:T5-5} shows the given and chosen APKs of the experiment with contrived data. \\
Generally, the remaining errors are situated in regions where we expect them to be. In the case of the approximation of surface data, \cref{fig:T4}, and the downward continuation of regularly distributed satellite data, \cref{fig:T12}, we find deviations to the solution in particular in the Andean region, the Himalayas and the Pacific Ring of Fire. This is reasonable as the gravitational potential contains much more local structure there which gets lost due to the noise and the ill-posedness. Further, we find that, in the case of approximating surface data, i.\,e.\ using potential data which is not damped due to satellite height, the methods obtain much better relative RMSEs while still counting more data errors than in the case of the downward continuation, see \cref{tab:STD}. The latter is clear, since more local structures are visible on the surface and appear relatively larger with respect to the noisy data. Obviously, they also need more iterations in this case. Again, as the data contains more information in this case this behaviour can be expected. Nonetheless, these experiments show that the LIPMP algorithms can indeed be used as standalone approximation algorithms. \\
Similar results are obtained for the downward continuation of irregularly distributed satellite data,  \cref{fig:T7}. However, we notice that some additional errors occur here in comparison to the results of the regularly distributed data, \cref{fig:T12}. In particular, these are located mostly in areas where we used less data, see, e.\,g., the North Atlantic and the Indian Ocean. This points out that the LIPMP algorithms are able to distinguish smoother and rougher regions on its own and, thus, prevent local gaps to have a global influence.\\
At last, we consider the test with the contrived data. First of all, we note that the LROFMP algorithm is able to approximate this data as well, see \cref{tab:STD}. The most important results are that the SHs are obtained exactly, the APKs are either clustered around the solutions or have a very small coefficient. For the latter, see \cref{fig:T5-5}. Note that those few wrongly chosen APKs may be caused by the present noise. The SHs are easier to distinguish, most likely because of their orthogonality. Hence, we see that the LROFMP algorithm is able to distinguish global trends and local anomalies.

\section{Conclusion and Outlook}
The downward continuation of the gravitational potential from satellite data is important for many reasons such as monitoring the climate change. One approach for this is presented by the IPMP algorithms. They seek iteratively the minimizer of the Tikhonov functional and, in this way, obtain a weighted linear combination in dictionary elements as an approximation. Thus, an a-priori chosen, finite dictionary is to a certain extent a bias for this approximation. The novel LIPMP algorithms include an add-on such that an infinite dictionary can be used. Further, a finite dictionary can be learnt as well.\\
Our numerical results show that both the non-learning IPMP as well as the LIPMP algorithms yield good results. However, the LIPMP algorithms have additional advantages in terms of CPU-runtime, storage demand, sparsity and the consequences of the number of different types of trial functions in use. Hence, we suggest that the IPMP algorithms may be used if those aspects are not critical because these methods are easier to implement. Otherwise, we advise to include the add-on, i.\,e. use the LIPMP algorithms. In particular, we prefer the LRFMP algorithm as it was used here and in \cite{Schneider2020} because it has a lower runtime than the LROFMP algorithm and is easier to implement.\\
In the future, we are interested to increase the number of used data as well as apply the algorithms to other geoscientific tasks, e.\,g., from seismology. 

\section*{Declarations}

\paragraph{Funding} The research was funded by the German Research Foundation (DFG; Deutsche Forschungsgemeinschaft), projects MI 655/7-2 and MI 655/14-1 (Naomi Schneider) and the University of Siegen (Volker Michel).
\paragraph{Conflicts of interest/Competing interests} Not applicable.
\paragraph{Availability of data, material and code} Publicly available data of GRACE and the EGM2008 was used, see  \cite{Devaraju2017,Flechtneretal2014,Flechtneretal2014-2,GRACEdata,Pavlisetal2012,Schmidtetal2008,Tapleyetal2004,GRACEdata2}. The particular datasets generated  and  analysed  during  the  current  study  are  available  from  the authors on reasonable request. The code generated and used during the current study is also available from the authors on reasonable request.
\paragraph{Author's contributions} The research was carried out for the dissertation of Naomi Schneider which was supervised by Volker Michel.

\footnotesize
\bibliography{biblio}
\normalfont

\newpage
\appendix
\section{Supplementary plots from the numerical experiments}
We give plots of examples of trial functions as well as the EGM2008 and the deviation from the mean field of the GRACE data in May 2008. Further, we give the irregular data grid and plots of the absolute approximation error of the comparison tests.
\begin{figure}[h]
		\begin{subfigure}{\textwidth}
		\centering
			\includegraphics[width=.49\textwidth]{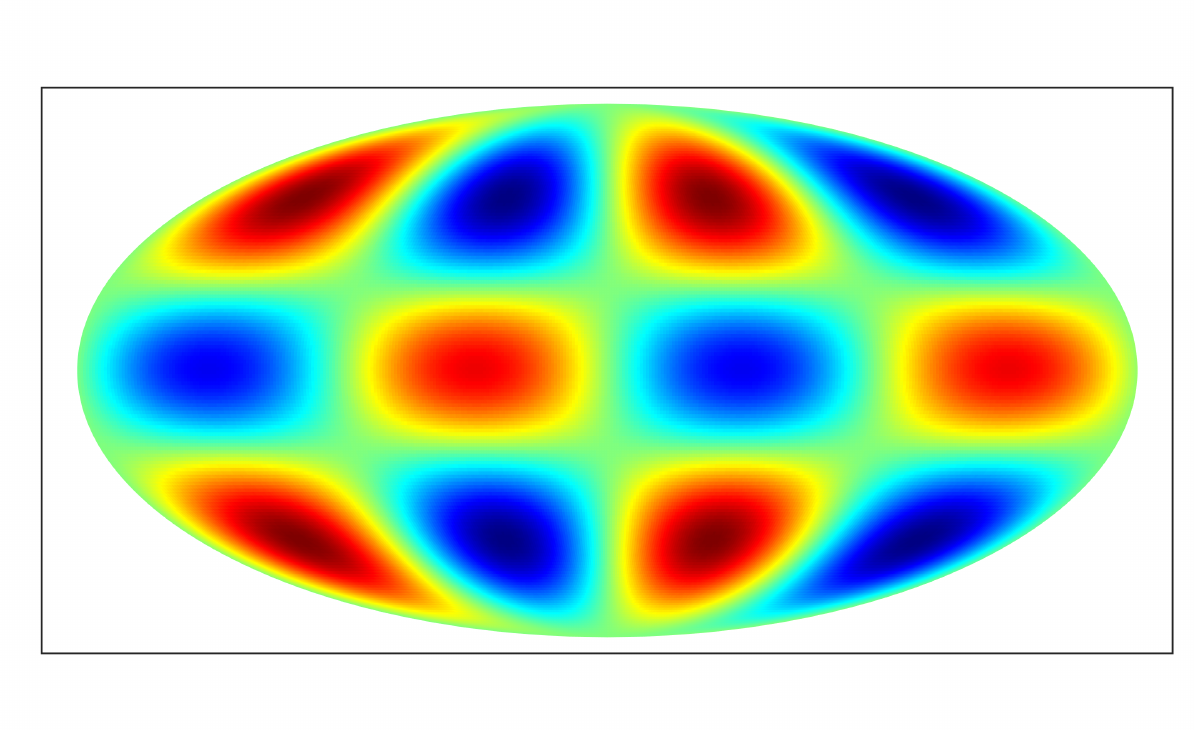}
			\includegraphics[width=.49\textwidth]{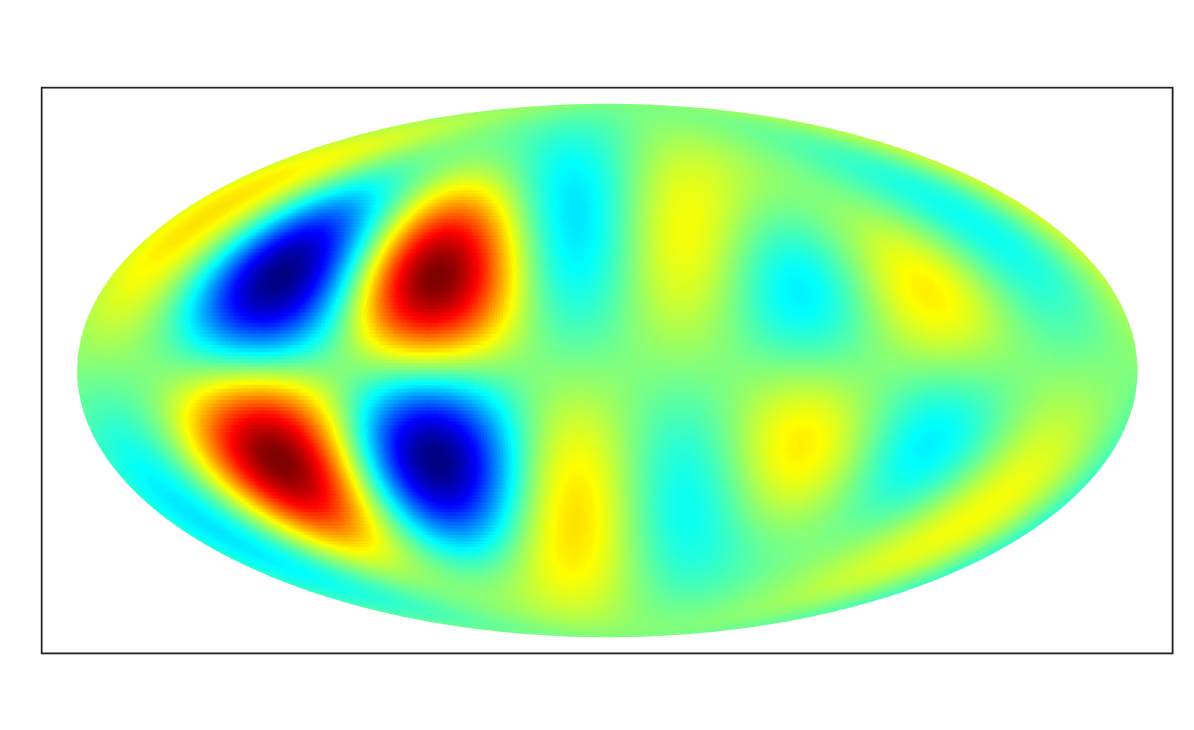}
			\caption{Example of a spherical harmonic (left) and a Slepian function (right).}
			\label{fig:TF1}
		\end{subfigure}
		\begin{subfigure}{\textwidth}
		\centering
			\includegraphics[width=.49\textwidth]{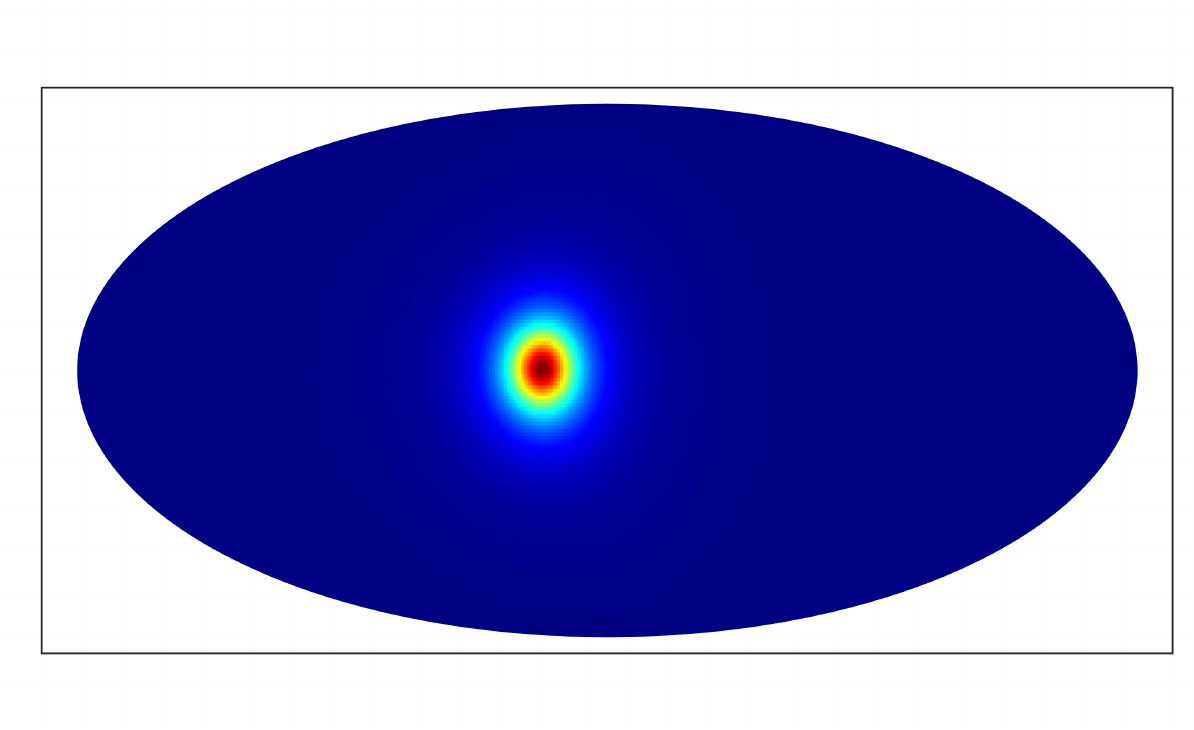}
			\includegraphics[width=.49\textwidth]{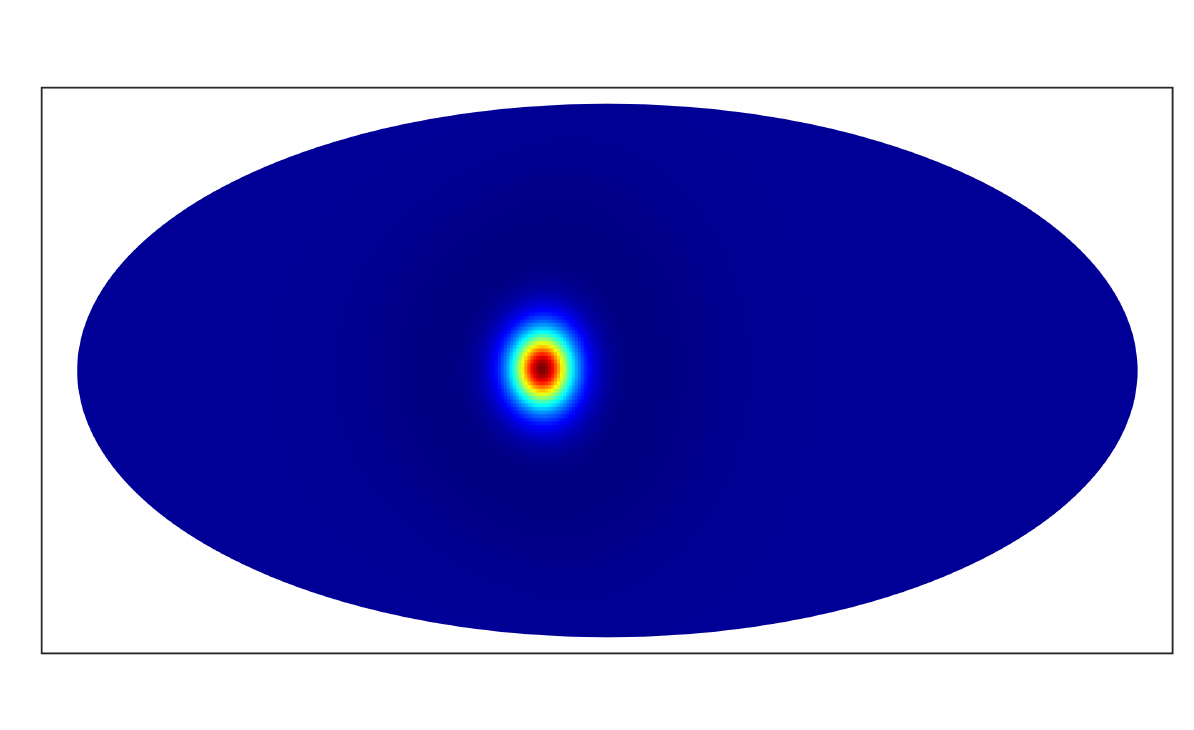}
			\caption{Example of an Abel--Poisson kernel (left) and wavelet (right).}
			\label{fig:TF2}
		\end{subfigure}
	\caption{Supplementary Plots: examples of the trial functions.}
	\label{fig:suppTF}
\end{figure}
\begin{figure}[h]
		\centering
			\includegraphics[width=.32\textwidth]{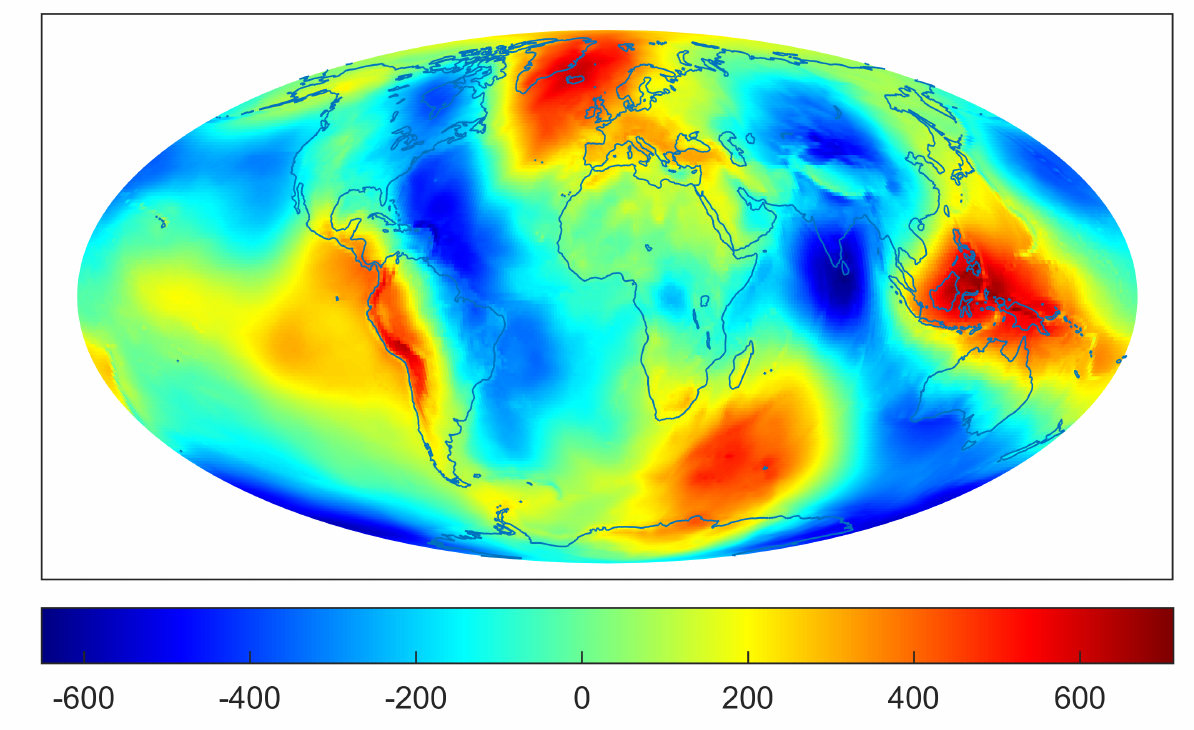}
			\includegraphics[width=.32\textwidth]{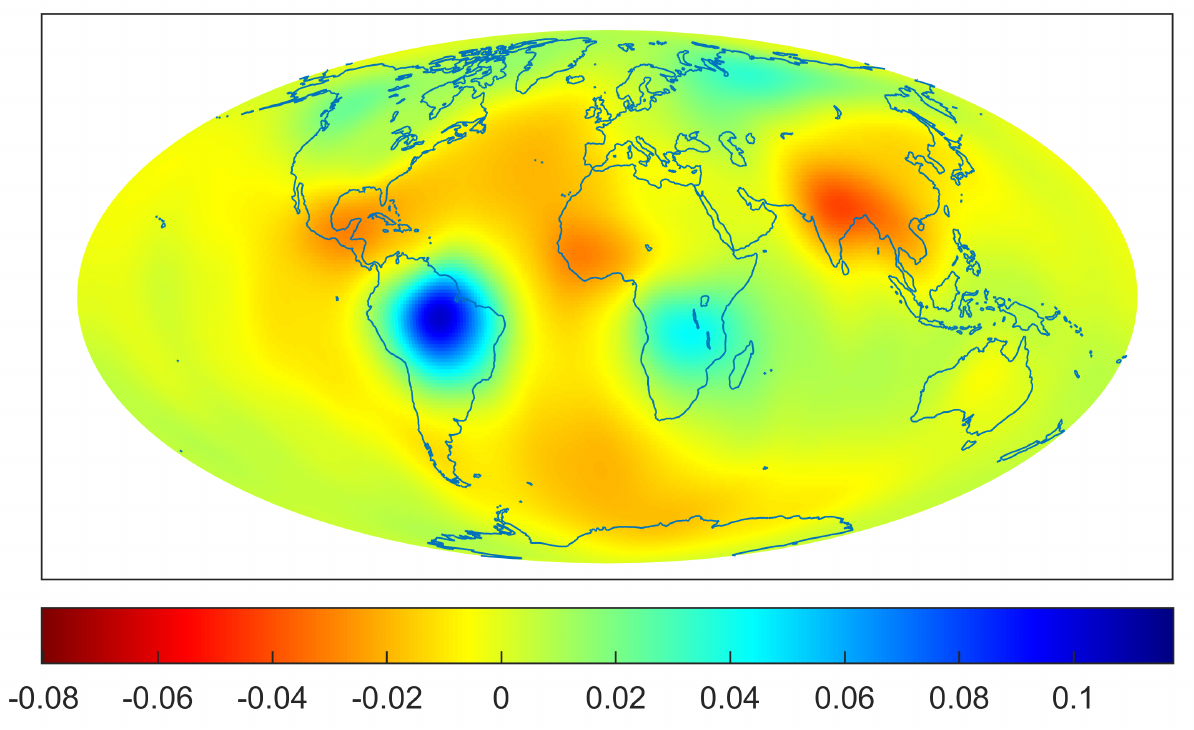}
		\includegraphics[width=.32\textwidth, scale=1.25]{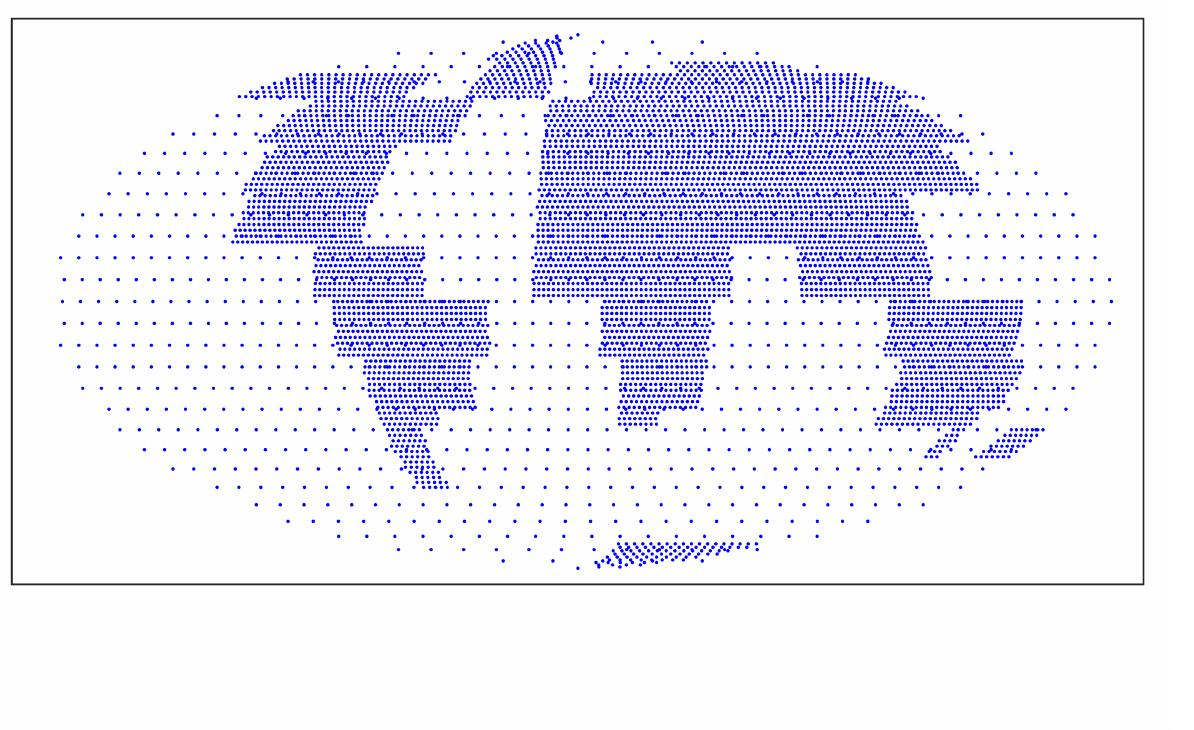}
			\caption{Supplementary plots: the data and the irregular data grid. EGM2008 (left; degree 3 to 2190) and deviation from the mean value of 2003 to 2013 in May 2008 due to GRACE (middle; degree 3 to 60). All values in $\mathrm{m}^2/\mathrm{s}^2$. Right: the irregular data grid.}
			\label{fig:dataandgrid}
\end{figure}
\begin{figure}
		\begin{subfigure}{\textwidth}
			\centering
			\includegraphics[width=.45\textwidth]{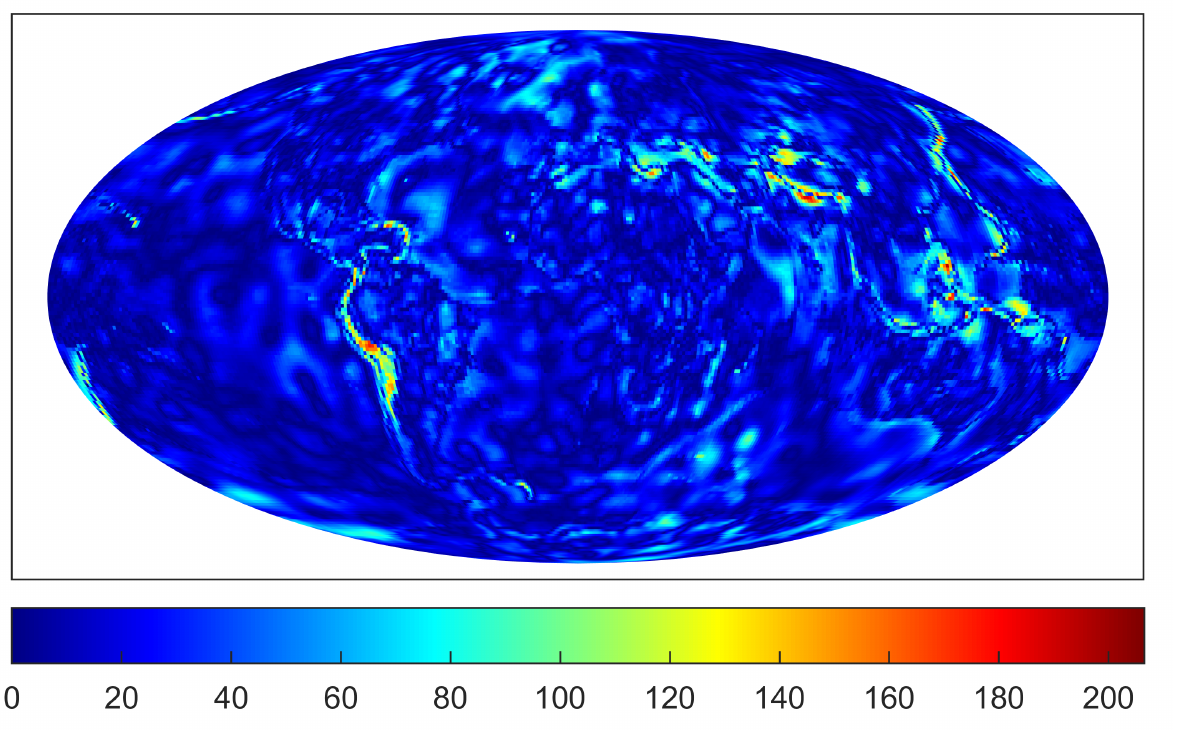}
			\includegraphics[width=.45\textwidth]{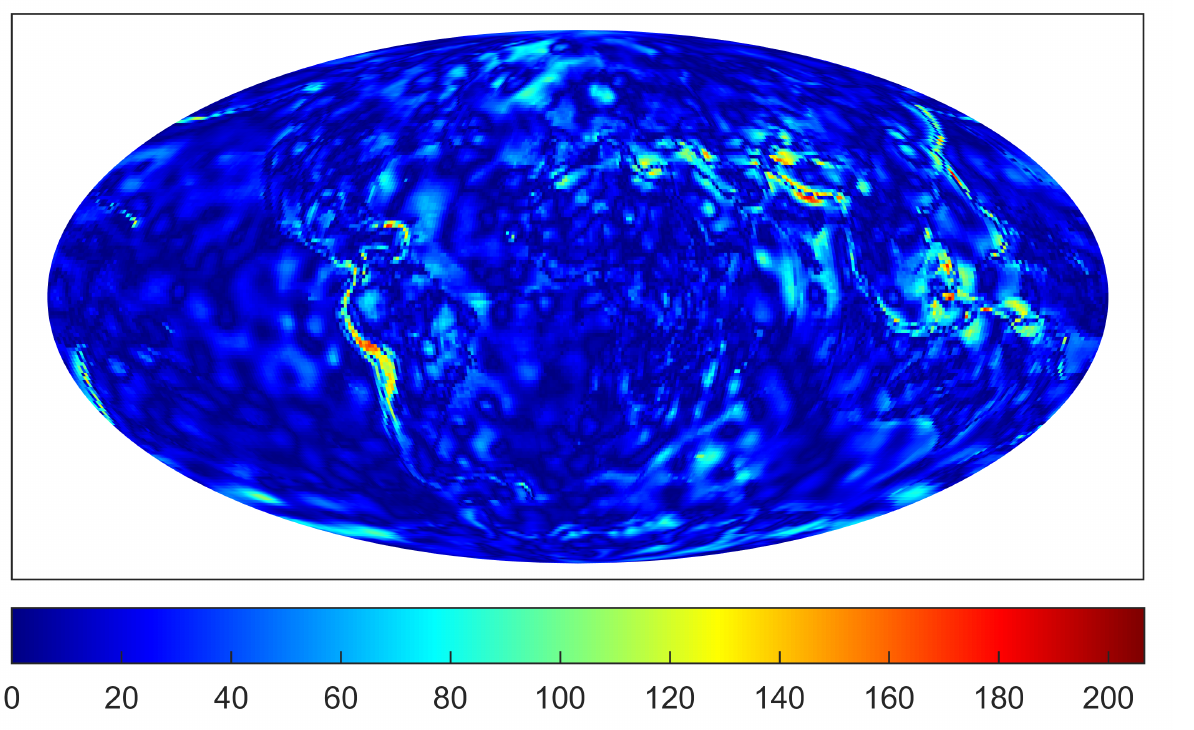}
			\includegraphics[width=.45\textwidth]{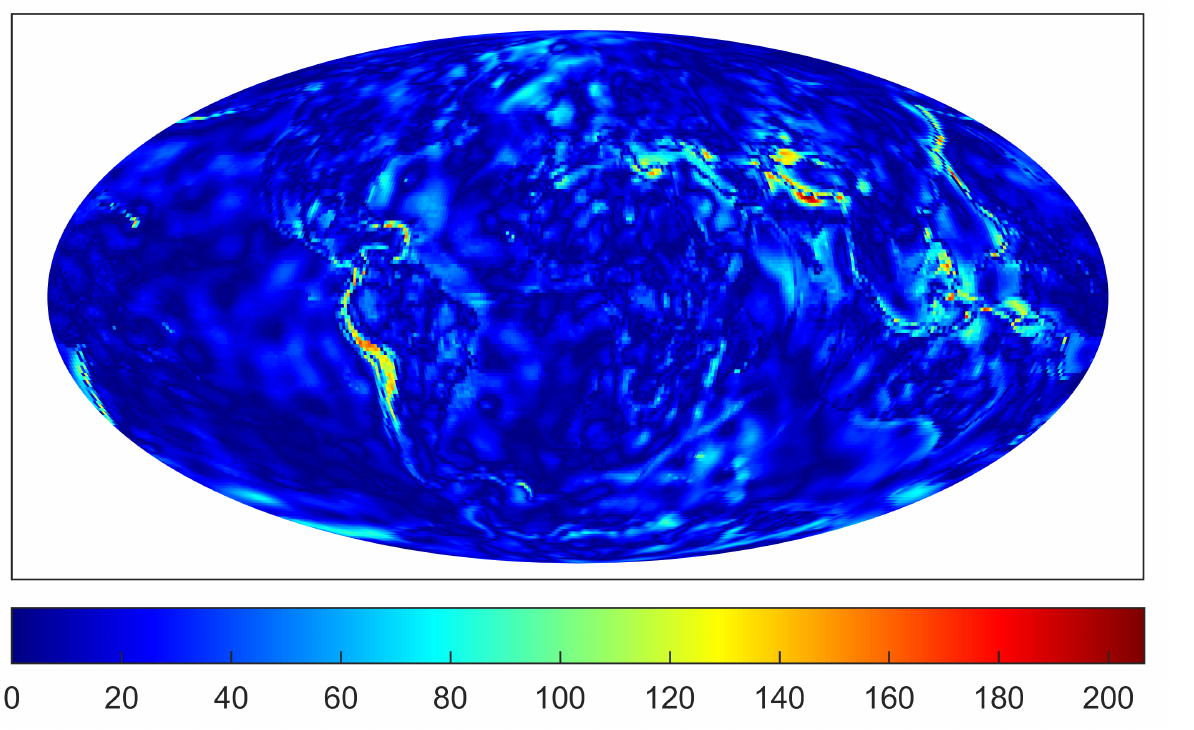}
			\includegraphics[width=.45\textwidth]{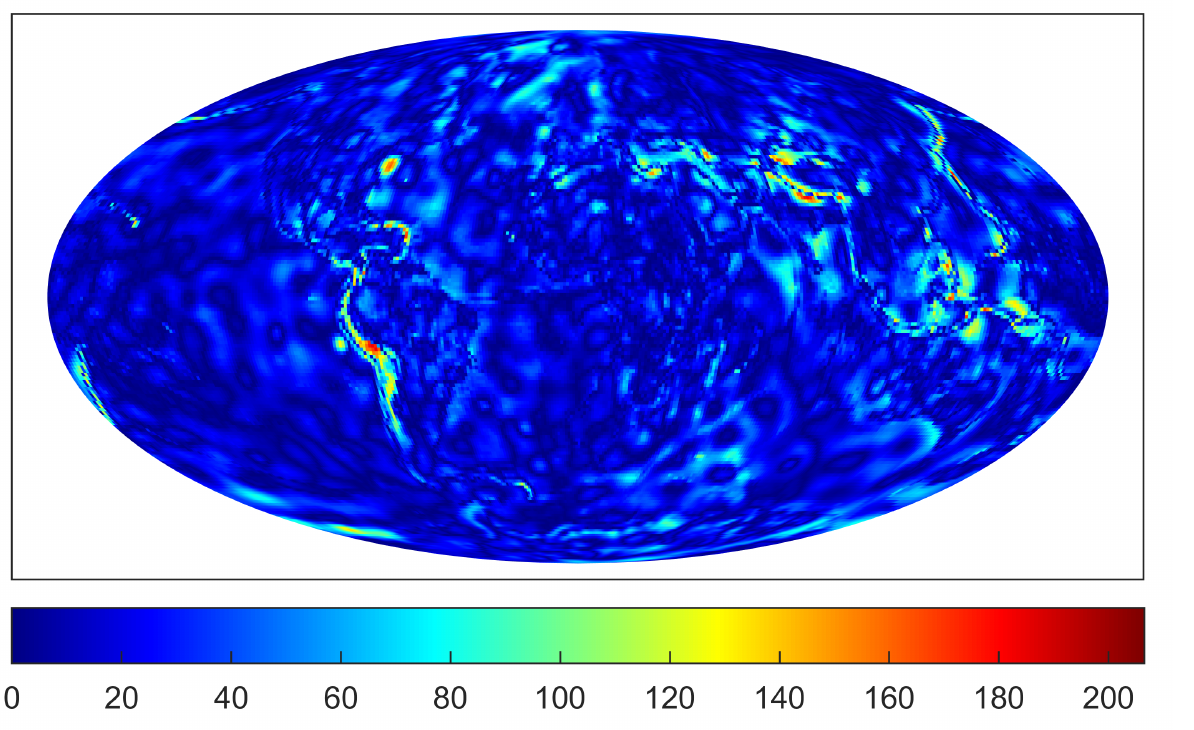}
			\caption{Absolute approximation errors obtained by the RFMP algorithm. The RFMP algorithm uses the manually chosen (left, upper row), the learnt (right, upper row), the non-stationary learnt (left, lower row) and the learnt-without-Slepian-functions (right, lower row) dictionary.}
			\label{fig:T1lrfmp:AErr}
		\end{subfigure}
		\begin{subfigure}{\textwidth}
			\centering
			\includegraphics[width=.45\textwidth]{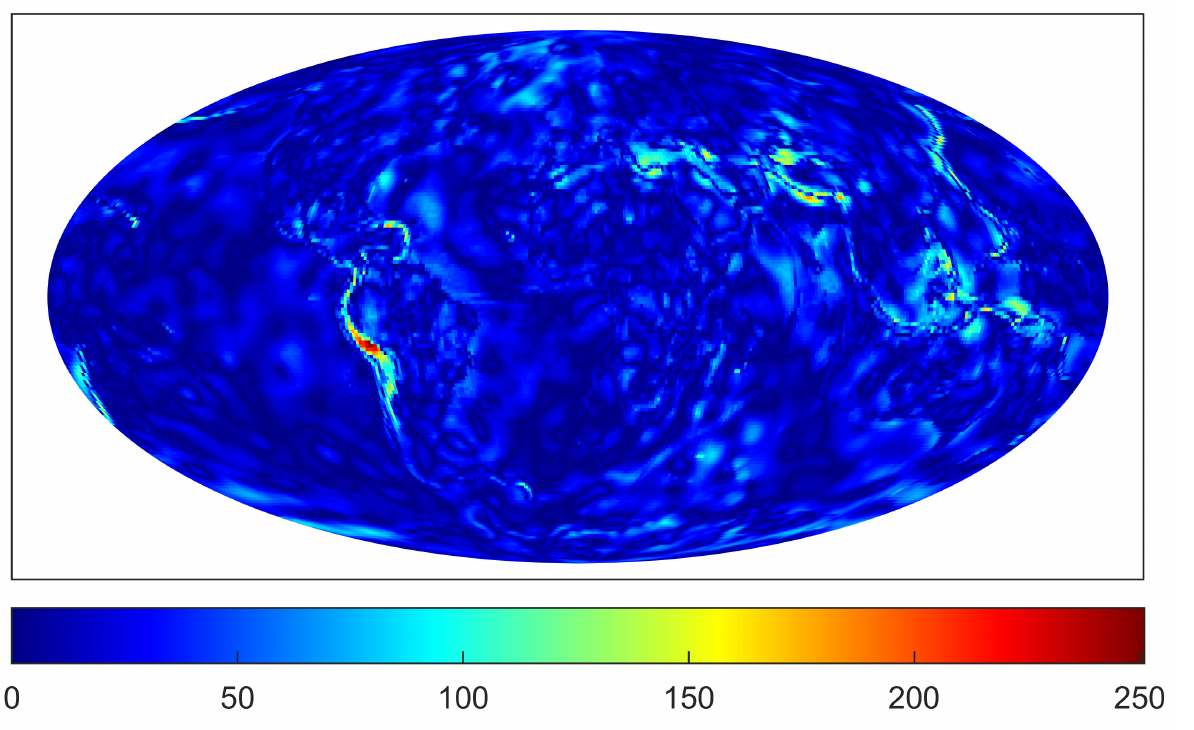}
			\includegraphics[width=.45\textwidth]{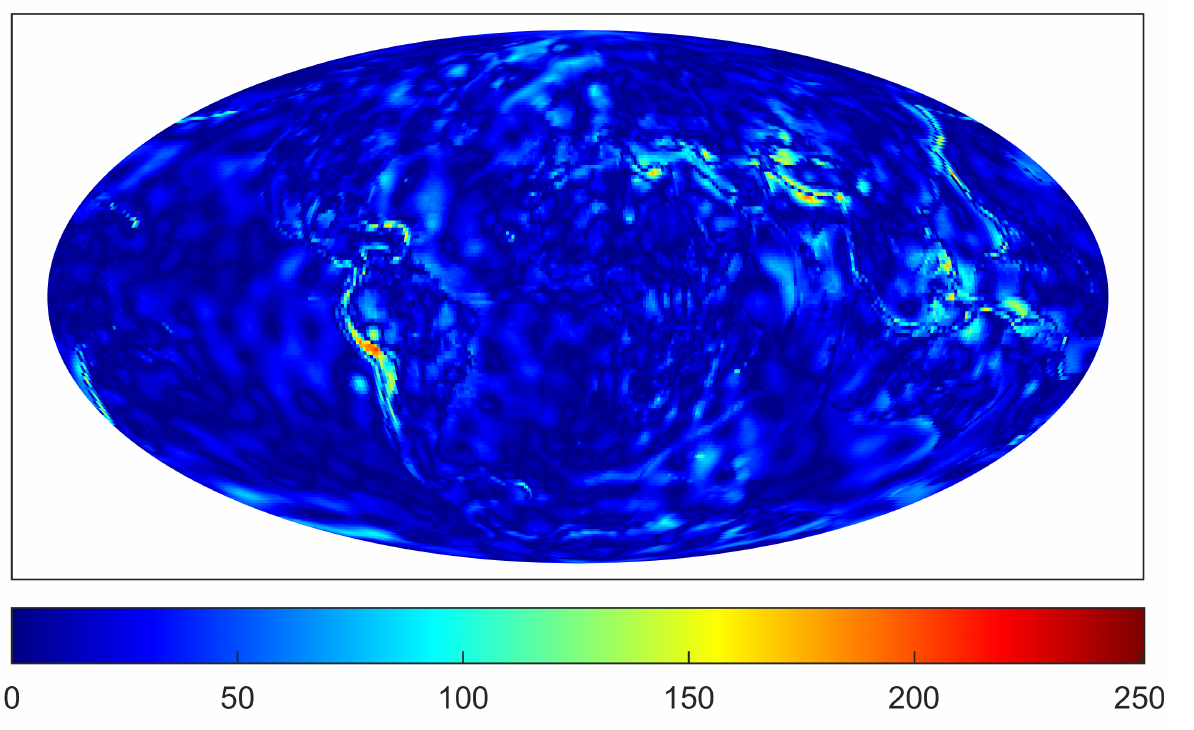}
			\includegraphics[width=.45\textwidth]{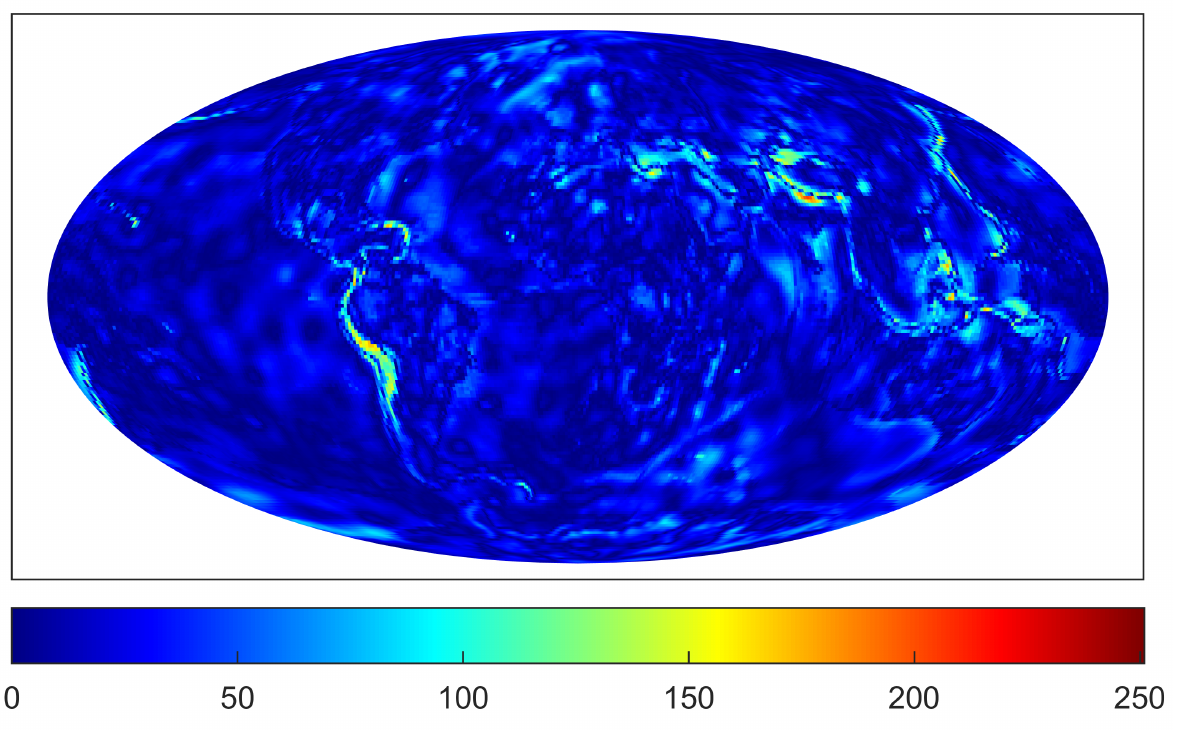}
			\includegraphics[width=.45\textwidth]{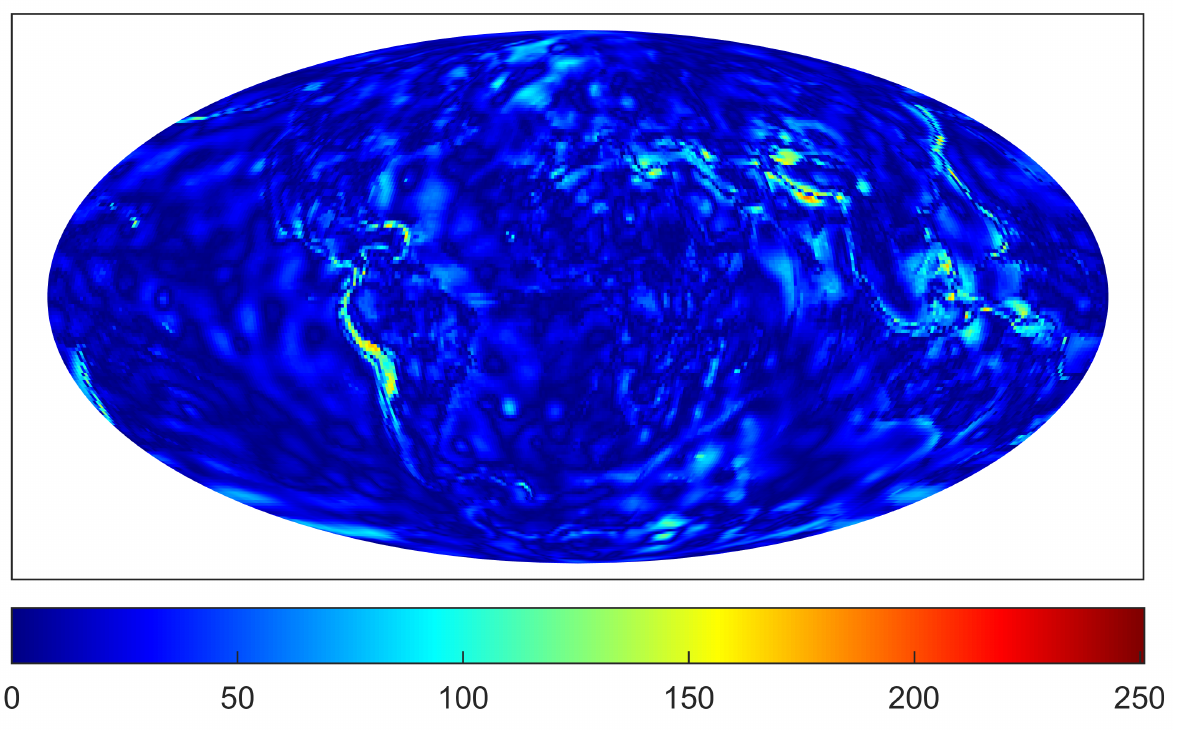}
			\caption{Absolute approximation errors obtained by the ROFMP algorithm. The ROFMP algorithm uses the manually chosen (left, upper row), the learnt (right, upper row), the non-stationary learnt (left, lower row) and the learnt-without-Slepian-functions (right, lower row) dictionary.}
			\label{fig:T1lrofmp:AErr}
		\end{subfigure}
\caption{Supplementary Plots: comparison of manually chosen and learnt dictionary with EGM2008 data. The colour scale is adapted for better comparability. All values in $\mathrm{m}^2/\mathrm{s}^2$. }
	\label{fig:supp2}
\end{figure}
\begin{figure}	
		\begin{subfigure}{\textwidth}
			\centering
			\includegraphics[width=.45\textwidth]{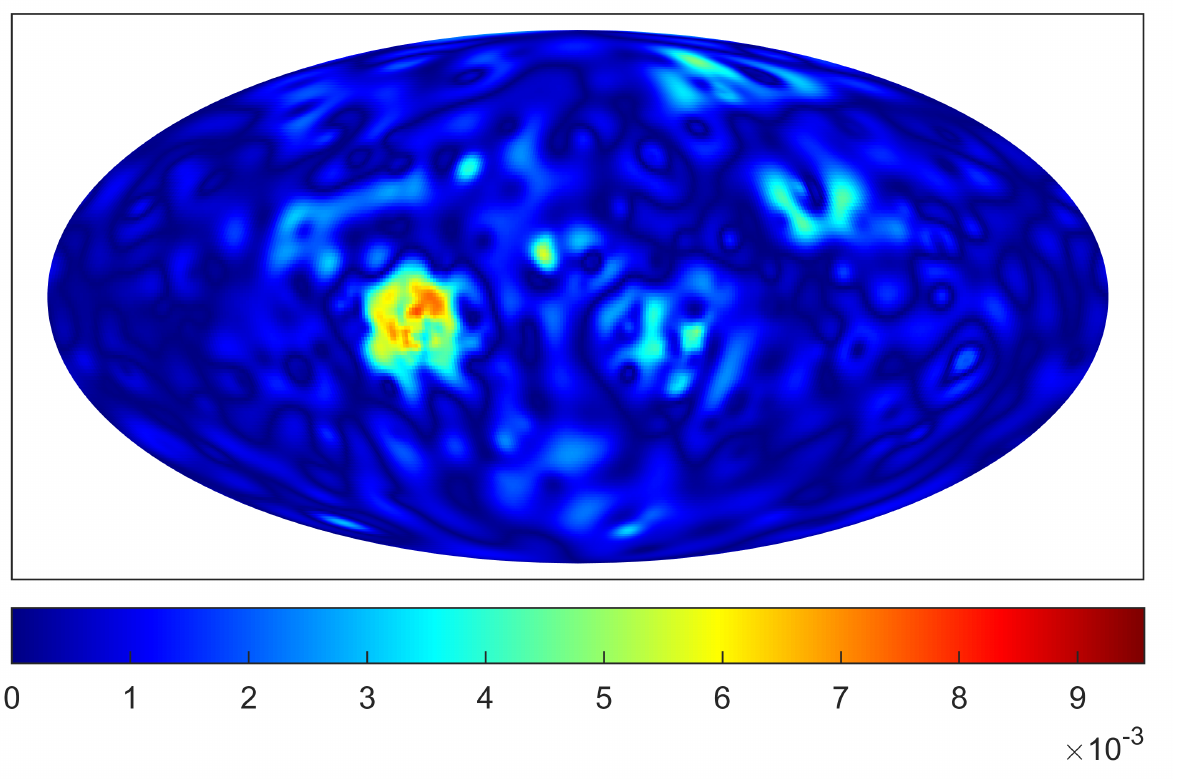}
			\includegraphics[width=.45\textwidth]{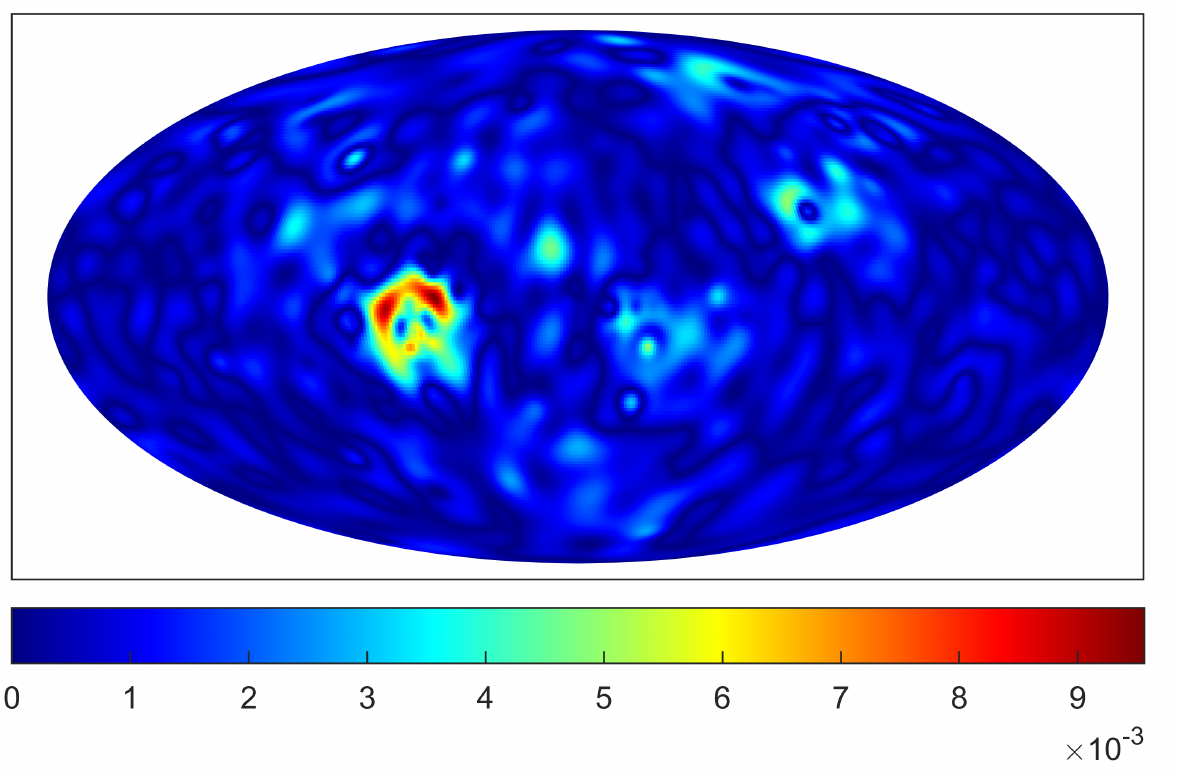}
			\includegraphics[width=.45\textwidth]{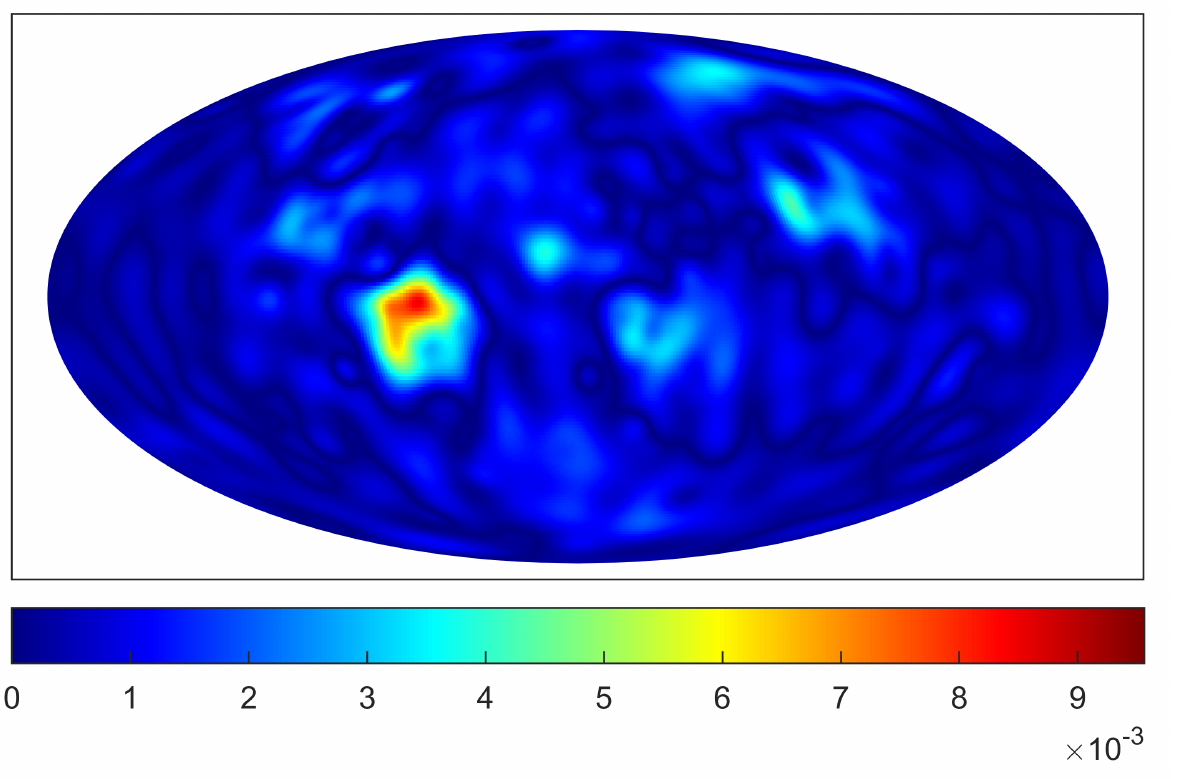}
			\includegraphics[width=.45\textwidth]{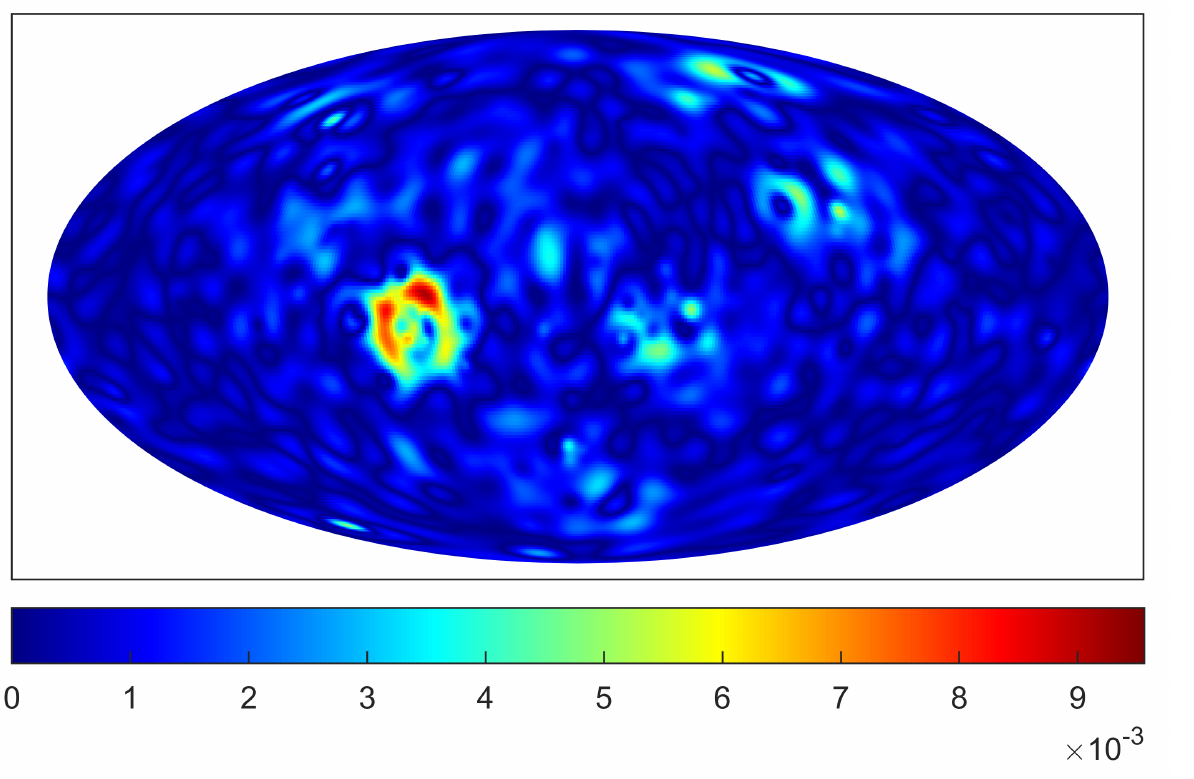}
			\caption{Absolute approximation error obtained by the RFMP algorithm. The RFMP algorithm uses the manually chosen (left, upper row), the learnt (right, upper row), the non-stationary learnt (left, lower row) and the learnt-without-Slepian-functions dictionary (right, lower row).}
			\label{fig:T2lrfmp:AErr}
		\end{subfigure}
		\begin{subfigure}{\textwidth}
			\centering
			\includegraphics[width=.45\textwidth]{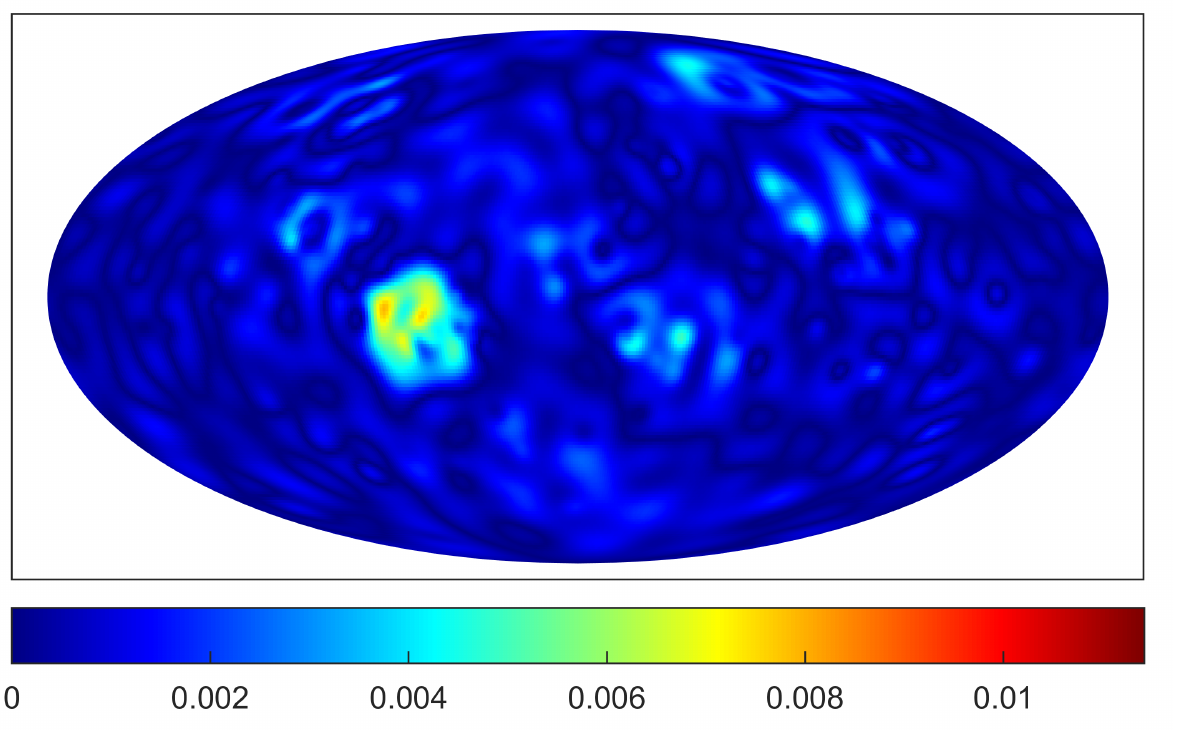}
			\includegraphics[width=.45\textwidth]{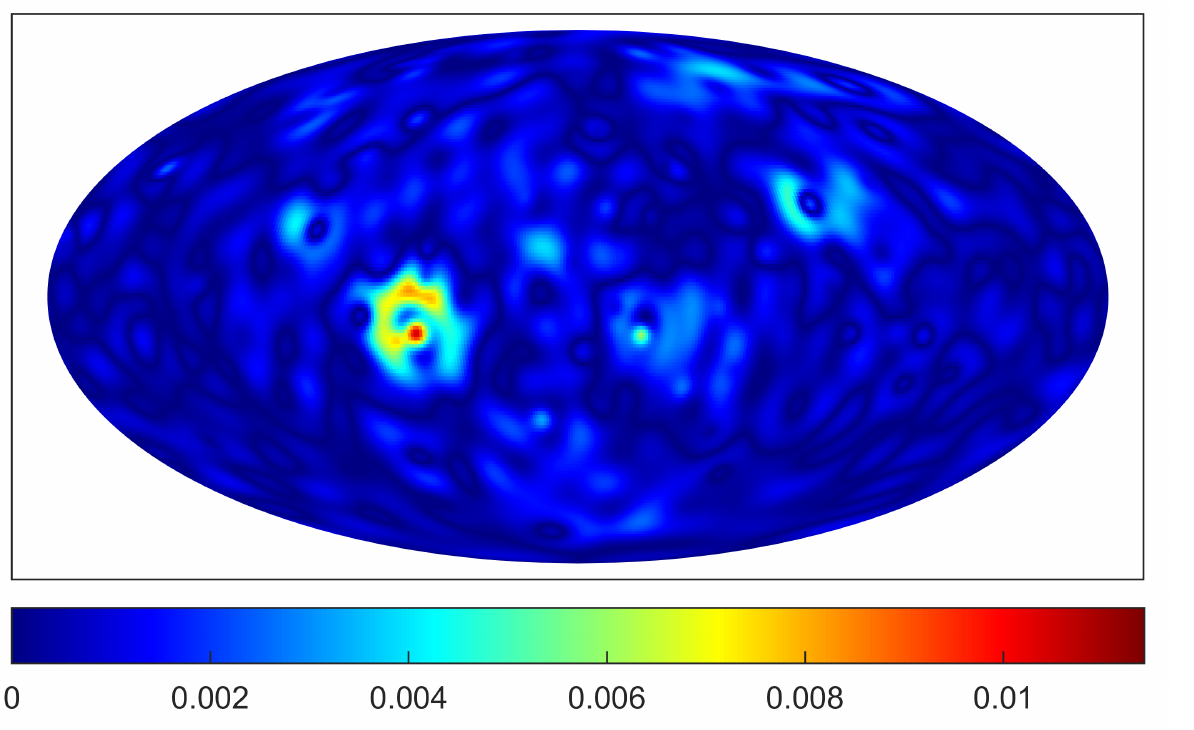}
			\includegraphics[width=.45\textwidth]{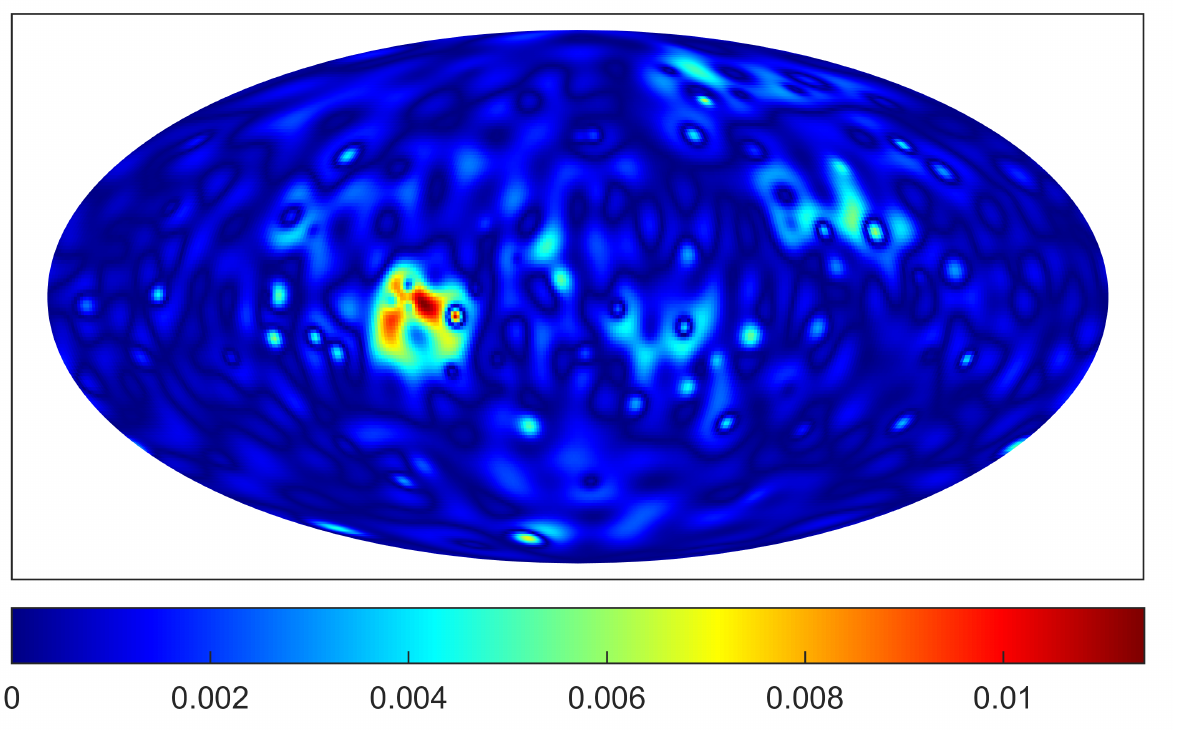}
			\includegraphics[width=.45\textwidth]{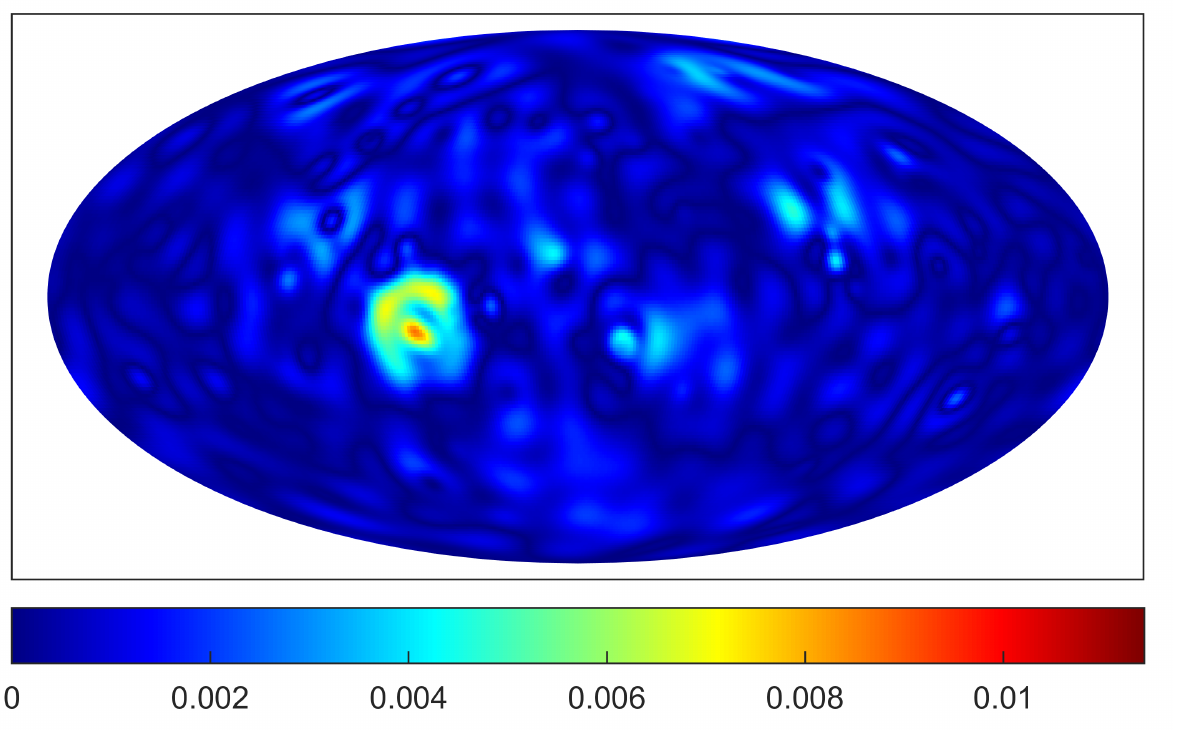}
			\caption{Absolute approximation error obtained by the ROFMP algorithm. The ROFMP algorithm uses the manually chosen (left, upper row), the learnt (right, upper row), the non-stationary learnt (left, lower row) and the learnt-without-Slepian-functions dictionary (right, lower row).}
			\label{fig:T2lrofmp:AErr}
		\end{subfigure}
	\caption{Supplementary Plots: comparison of manually chosen and learnt dictionary with GRACE data (May 2008). The colour scale is adapted for better comparability. All values in $\mathrm{m}^2/\mathrm{s}^2$.}
	\label{fig:supp3}
\end{figure}

\end{document}